\documentclass{amsart}%
\RequirePackage{ifthen}
\newboolean{omarfont}
\setboolean{omarfont}{false}
\newboolean{usesvjour}
\setboolean{usesvjour}{false}%
    \RequirePackage{ifthen}
    \provideboolean{useutopia}
    \setboolean{useutopia}{false}%
    \provideboolean{usesvjour}
    \setboolean{usesvjour}{false}%
    \ifthenelse{\boolean{useutopia}}{%
      \RequirePackage{ifthen}
      \RequirePackage{iftex}
      \provideboolean{usemathrsfs}%
      \provideboolean{useutopia}
      \setboolean{useutopia}{false}%
      \ifXeTeX
        \RequirePackage[american,british]{babel}
        \RequirePackage{xltxtra}%
        \DeclareSymbolFont{usualmathcal}{OMS}{cmsy}{m}{n}
        \DeclareSymbolFontAlphabet{\mathcalbf}{usualmathcal}
        \ifthenelse{\boolean{useutopia}}{
          \RequirePackage[utopia,euro]{mathdesign}
          \RequirePackage[OMLmathbf,OMLmathsfit]{isomath}
        }{
          \RequirePackage{amssymb}    %
          \RequirePackage{bbold}    %
        }
        \RequirePackage{mathrsfs} %
        \setboolean{usemathrsfs}{true}%
      \else%
        \RequirePackage[utf8]{inputenc}%
        \ifthenelse{\boolean{useutopia}}{%
          \RequirePackage[utopia,euro]{mathdesign}%
          \RequirePackage[OMLmathrm,OMLmathbf,OMLmathsfit]{isomath}
          \DeclareSymbolFont{usualmathcal}{OMS}{cmsy}{m}{n}
          \DeclareSymbolFontAlphabet{\mathcalbf}{usualmathcal}
          \RequirePackage{stmaryrd} %
          \providecommand{\diracdelta}[1][]{\ensuremath{\deltaup_{#1}}}
          
          \providecommand{\lap}{\ensuremath{\Deltaup}}
          \providecommand{\measure}[1]{\ensuremath{\mathcalbf{\uppercase{#1}}}}
        }{%
          \RequirePackage{amssymb}    %
          \RequirePackage{mathrsfs} %
          \RequirePackage{bbold}    %
          \RequirePackage{stmaryrd} %
          \setboolean{usemathrsfs}{true}%
          \providecommand{\mathcalbf}{\mathcal}
        }%
        \RequirePackage[american,british]{babel}
      \fi%
      \RequirePackage{mathtools}%
      \RequirePackage{nicefrac}
      \RequirePackage{stmaryrd} %
      \RequirePackage{amsthm}%
      \RequirePackage{enumerate}%
      \RequirePackage{xspace}%
      \RequirePackage{verbatim}%
      \RequirePackage{fancyvrb}
      \RequirePackage{hyperref}%
      \RequirePackage{listings}%
      \RequirePackage{xcolor}
      \RequirePackage{epsfig}
      \RequirePackage{graphicx}
      \RequirePackage{booktabs}%
      \RequirePackage{tikz}%
      \usetikzlibrary{calc}%
      \provideboolean{usebeamer}%
      \provideboolean{isthesis}%
      \provideboolean{isamsltex}%
      \provideboolean{issiamltex}%
    }{%
      \RequirePackage{ifthen}
      \RequirePackage{iftex}
      \provideboolean{usemathrsfs}%
      \provideboolean{useutopia}
      \setboolean{useutopia}{false}%
      \ifXeTeX
        \RequirePackage[american,british]{babel}
        \RequirePackage{xltxtra}%
        \DeclareSymbolFont{usualmathcal}{OMS}{cmsy}{m}{n}
        \DeclareSymbolFontAlphabet{\mathcalbf}{usualmathcal}
        \ifthenelse{\boolean{useutopia}}{
          \RequirePackage[utopia,euro]{mathdesign}
          \RequirePackage[OMLmathbf,OMLmathsfit]{isomath}
        }{
          \RequirePackage{amssymb}    %
          \RequirePackage{bbold}    %
        }
        \RequirePackage{mathrsfs} %
        \setboolean{usemathrsfs}{true}%
      \else%
        \RequirePackage[utf8]{inputenc}%
        \ifthenelse{\boolean{useutopia}}{%
          \RequirePackage[utopia,euro]{mathdesign}%
          \RequirePackage[OMLmathrm,OMLmathbf,OMLmathsfit]{isomath}
          \DeclareSymbolFont{usualmathcal}{OMS}{cmsy}{m}{n}
          \DeclareSymbolFontAlphabet{\mathcalbf}{usualmathcal}
          \RequirePackage{stmaryrd} %
          \providecommand{\diracdelta}[1][]{\ensuremath{\deltaup_{#1}}}
          
          \providecommand{\lap}{\ensuremath{\Deltaup}}
          \providecommand{\measure}[1]{\ensuremath{\mathcalbf{\uppercase{#1}}}}
        }{%
          \RequirePackage{amssymb}    %
          \RequirePackage{mathrsfs} %
          \RequirePackage{bbold}    %
          \RequirePackage{stmaryrd} %
          \setboolean{usemathrsfs}{true}%
          \providecommand{\mathcalbf}{\mathcal}
        }%
        \RequirePackage[american,british]{babel}
      \fi%
      \RequirePackage{mathtools}%
      \RequirePackage{nicefrac}
      \RequirePackage{stmaryrd} %
      \RequirePackage{amsthm}%
      \RequirePackage{enumerate}%
      \RequirePackage{xspace}%
      \RequirePackage{verbatim}%
      \RequirePackage{fancyvrb}
      \RequirePackage{hyperref}%
      \RequirePackage{listings}%
      \RequirePackage{xcolor}
      \RequirePackage{epsfig}
      \RequirePackage{graphicx}
      \RequirePackage{booktabs}%
      \RequirePackage{tikz}%
      \usetikzlibrary{calc}%
      \provideboolean{usebeamer}%
      \provideboolean{isthesis}%
      \provideboolean{isamsltex}%
      \provideboolean{issiamltex}%
    }
    \colorlet{a}{magenta}
    \colorlet{b}{green!75!blue}
    \colorlet{c}{yellow!87.5!red}
    \colorlet{d}{cyan}
    \colorlet{e}{red}
    \colorlet{f}{blue}
    \colorlet{g}{white}
    \colorlet{h}{black!50}
    \colorlet{i}{black}
    \colorlet{j}{black!75}

    \providecommand{\linkedurl}[1]{\url{#1}}
    \providecommand{\linkedemail}[1]{\href{mailto:#1}{#1}}
    
    \providecommand{\email}[1]{{\linkedemail{#1}}}

    \providecommand{\Ignore}[1]{}
    \providecommand{\ignore}[1]{}
    \providecommand{\freeze}[1]{}%
    \providecommand{\crossout}[1]{{\textcolor{i!20}{#1}}}
    \providecommand{\highlightcolor}{a}
    \providecommand{\highlight}[1]{{\color{\highlightcolor}#1}}

    \providecommand{\memo}[1]{%
      \ensuremath{%
        \framebox{\tiny\textbf{\kern-2pt\textsf{#1}}\kern-2pt}%
      }%
      \xspace}

    \RequirePackage{alphalph}
    \provideboolean{shownotes}
    \setboolean{shownotes}{true}%
    \newcounter{margnote}[page]
    \providecommand{\mgcolor}{a}
    \providecommand{\mgcolorset}[1]{\renewcommand{\mgcolor}{\alphalph{#1}}}
    \providecommand{\mgcolorsetbycounter}[1]{%
      \ifthenelse{\value{#1}<11}{%
        \renewcommand{\mgcolor}{\alph{#1}}%
      }{%
        \renewcommand{\mgcolor}{a}}%
    }
    \providecommand{\mgcolormake}{\mgcolorsetbycounter{margnote}}
    \providecommand{\mgcolorstepby}[1]{
      \setcounter{tmpcounter}{\value{margnote}}%
      \addtocounter{tmpcounter}{#1}%
      \mgcolorsetbycounter{tmpcounter}%
    }%
    \providecommand{\margnotecolor}{%
      \ifthenelse{\value{margnote}=0}{%
        \mgcolorset{10}
      }{%
        \ifthenelse{\value{margnote}<7}{%
          \mgcolormake%
        }{%
          \ifthenelse{\value{margnote}=7}{\mgcolorset{10}}{%
            \ifthenelse{\value{margnote}<11}{\mgcolormake}{%
              \ifthenelse{\value{margnote}<17}{\mgcolorstepby{-10}}{%
                \mgcolorset{10}%
              }%
            }%
          }%
        }%
      }%
    }%
    \providecommand{\margnotemark}{{\colorbox{\mgcolor}{\tiny\color{g}\upshape\texttt{\arabic{page}.\arabic{margnote}}}}}
    \providecommand{\margnote}[2][]{%
      \ifthenelse{%
        \boolean{shownotes}%
      }{%
        \stepcounter{margnote}%
        \margnotecolor%
        \margnotemark %
        \marginpar{%
          \color{\mgcolor}%
          \texttt{%
            \begin{minipage}{2cm}%
              \raggedright\tiny%
              \margnotemark%
              #2%
              \\
              {\ifx|#1|{}\else{ - #1}\fi}%
            \end{minipage}%
          }%
        }%
      }{%
      }%
    }%
    \providecommand{\mathnote}[2][]{%
      \ifthenelse{%
        \boolean{shownotes}%
      }{%
        \stepcounter{margnote}%
        \margnotecolor%
        \text{%
          \colorbox{\mgcolor}{%
            \color{g}%
            \texttt{%
              \tiny%
                  \margnotemark: %
                  \ifx|#1|{}\else{#1:}\fi%
                  #2%
            }%
          }%
        }%
      }{%
      }%
    }%
    \providecommand{\textnote}[2][]{%
      \ifthenelse{%
        \boolean{shownotes}%
      }{%
        \stepcounter{margnote}%
        \margnotecolor%
        \ \\
        \text{%
          \colorbox{\mgcolor}{%
            \begin{minipage}{.9\textwidth}
            \color{g}%
            \texttt{%
              \margnotemark: %
              \ifx|#1|{}\else{#1: }\fi%
              #2%
            }%
            \end{minipage}
          }%
        }%
      }{%
      }%
    }%

    \providecommand{\todo}[2][]{\margnote[#1]{To do: #2}}
    \providecommand{\Todo}[1]{
      \ifthenelse{\boolean{shownotes}}{
        \begin{center}
        \begin{tikzpicture}
         \node[fill=a!17]{
           \begin{minipage}{\textwidth}
             \texttt{To do:}
             \\
             \texttt{\bfseries{\small #1}}
           \end{minipage}
         };
        \end{tikzpicture}
        \end{center}
      }{}}
    \provideboolean{showrevisions}
    \setboolean{showrevisions}{true}
    \newcommand{\revisionsheader}{***\newline\Warning{the following part is under revision}}
    \newcommand{\revisionsfooter}{***\newline\Warning{end of part under revision}}

    \providecommand{\Warning}[1]{    
      \begin{tikzpicture}
        \node[fill=a!27]{
          \begin{minipage}{\textwidth}
            \texttt{\bfseries{\small Warning: #1}}
          \end{minipage}
        };
      \end{tikzpicture}
    }

    \provideboolean{showcomments}
    \providecommand{\margincomment}[1]{
    \ifthenelse{\boolean{showcomments}}{\marginpar{\tiny #1}}{}
    }
    \provideboolean{showchanges}
    \setboolean{showchanges}{false}
    \providecommand{\changes}[2][]{%
      \ifthenelse{\boolean{showchanges}}{{\ifx|#1|{}\else\margnote{#1}\fi\highlight{#2}}}{#2}}
    \providecommand{\mathchanges}[2][]{%
      \ifthenelse{\boolean{showchanges}}{{\ifx|#1|{}\else\mathnote{#1}\fi\highlight{#2}}}{#2}}

    \providecommand{\changefromto}[3][replace with]{%
      \ifthenelse{\boolean{showchanges}}{%
        {\crossout{#2}\margnote{#1}}{\highlight{#3}}}{%
        #3\xspace}%
    }
    \providecommand{\ChangePar}[3][]{%
      \ifthenelse{\boolean{showchanges}}{
        {\par\textcolor{i!20}{#2}\ifx|#1|\else\margnote{#1}\fi}{\par\textcolor{a}{#3}}
      }{%
        \par #3%
      }%
    }
    \providecommand{\InsertPar}[1]{
      \ifthenelse{\boolean{showchanges}}
      {{\par$\mapsto$ \textcolor{blue}{#1}}}
      {\par #1}
    }
    
    \providecommand{\mathchangefromto}[3][]{\crossout{#2}\ifx|#1|\else\mathnote{#1}\fi\highlight{#3}}

    \providecommand{\mathscript}
    	   {\mathscr}

     \providecommand{\cA}{\ensuremath{\mathscript A}\xspace}

     \providecommand{\bbbold}{\mathbb}

     \providecommand{\rN}{\ensuremath{\bbbold N}\xspace}
     
     \providecommand{\rP}{\ensuremath{\bbbold P}\xspace}
     
     \providecommand{\rR}{\ensuremath{\bbbold R}\xspace}
     
     \providecommand{\rT}{\ensuremath{\bbbold T}\xspace}

    \providecommand{\Ae}[1][]{\ensuremath{\ifx|#1|{\ }\else{\:#1\text{-}}\fi\text{almost everywhere }}\xspace}
    \providecommand{\Aa}[1][]{\ensuremath{\text{ for }\ifx|#1|{}\else{\:#1\text{-}}\fi\text{almost all }}}
    \providecommand{\as}[1][]{\ensuremath{\ifx|#1|{\ }\else{#1\text{-}}\fi\text{almost surely}}\xspace}
    \providecommand{\aposteriori}{aposteriori\xspace}
    
    \providecommand{\apriori}{{apriori}\xspace}

     \providecommand{\naturals}{\rN\xspace}

     \providecommand{\reals}{\rR}

     \providecommand{\R}[1]{\reals^{#1}}
     
     \providecommand{\fieldmats}[3][F]{\csname#1\endcsname{#2\times#3}}
     \providecommand{\realmats}[2]{\fieldmats[R]{#1}{#2}}

     \providecommand{\torus}[1]{\rT\ifthenelse{\equal{#1}1}{}{^#1}}

     \providecommand{\diracdelta}[1][]{\ensuremath{{\mathrm{\delta}}\ifx|#1|{}\else_{#1}\fi}}

     \providecommand{\pic}{\ensuremath{\mathrm\pi}}
     \providecommand{\pifracl}[2][]{\fracl{\ifx|#1|\else#1\fi\pic}{#2}}
     \providecommand{\pifrac}[2][]{\frac{\ifx|#1|\else#1\fi\pic}{#2}}

     \providecommand{\inner}{\cdot}
     \providecommand{\outerp}{\wedge}

     \providecommand{\frobinner}{:}

     \providecommand{\W}{\ensuremath{\varOmega}\xspace}

     \providecommand{\w}{\ensuremath{\omega}\xspace}

     \providecommand{\qp}[2][]{\ifx|#1|\left(\else\csname#1\endcsname(\fi{#2}\ifx|#1|\right)\else\csname#1\endcsname)\fi}

     \providecommand{\qpreg}[1]{\ensuremath{(#1)}}
     \providecommand{\qpbig}[1]{\qp[big]{#1}}%
     \providecommand{\qpBig}[1]{\ensuremath{\Big(#1\Big)}}
     \providecommand{\qpbigg}[1]{\ensuremath{\bigg(\!#1\!\bigg)}}
     \providecommand{\qpBigg}[1]{\ensuremath{\Bigg(\!#1\!\Bigg)}}
     \providecommand{\qb}[2][]{\ifx|#1|\left[\else\csname#1\endcsname[\fi{#2}\ifx|#1|\right]\else\csname#1\endcsname]\fi}
     \providecommand{\qc}[2][]{\ifx|#1|\left\{\else\csname#1\endcsname\{\fi{#2}\ifx|#1|\right\}\else\csname#1\endcsname\}\fi}

     \providecommand{\qa}[1]{\ensuremath{\left\langle{#1}\right\rangle}}
     \providecommand{\qareg}[1]{\ensuremath{\langle#1\rangle}}
     \providecommand{\qabig}[1]{\ensuremath{\big\langle#1\big\rangle}}
     \providecommand{\qaBig}[1]{\ensuremath{\Big\langle#1\Big\rangle}}
     \providecommand{\qabigg}[1]{\ensuremath{\bigg\langle#1\bigg\rangle}}
     \providecommand{\qaBigg}[1]{\ensuremath{\Bigg\langle#1\Bigg\rangle}}
     \providecommand{\opinter}[2]{\ensuremath{\left(#1,#2\right)}\xspace}
     \providecommand{\clinter}[2]{\ensuremath{\left[#1,#2\right]}\xspace}

     \providecommand{\compowqp}[2]{\ensuremath{\qp{\!#2\!\!}^{\kern -.4em #1}\!}}
     
     \providecommand{\powqpreg}[2]{\ensuremath{%
         \qpreg{#2}^{\kern 0em\lower .1ex\hbox{\scriptsize $#1$}}\kern-.3em}}
     \providecommand{\powqpbig}[2]{\ensuremath{%
         \qpbig{#2}^{\kern -.2em\lower .3ex\hbox{\scriptsize $#1$}}\kern-.3em}}
     \providecommand{\powqpBig}[2]{\ensuremath{%
         \qpBig{#2}^{\kern -.2em\lower .3ex\hbox{\scriptsize $#1$}}\kern-.3em}}
     \providecommand{\powqpbigg}[2]{\ensuremath{%
         \qpbigg{#2}^{\kern -.2em\lower .3ex\hbox{\scriptsize $#1$}}\kern-.3em}}
     \providecommand{\powqpBigg}[2]{\ensuremath{%
         \qpBigg{#2}^{\kern -.2em\lower .3ex\hbox{\scriptsize $#1$}}}}
     \providecommand{\powp}[3][]{#3\ifx|#1|^{#2}\else{#1}^{#2}\fi}%
     \providecommand{\pow}[2][]{\ifx|#1|\operatorname{pow}^{#2}\else\powp{#2}{#1}\fi}%
     \providecommand{\ppow}[3][]{\powp[#1]{#3}{#2}}

     \providecommand{\powsqrt}[2][2]{\powp{\fracl1{#1}}{#2}}

     \providecommand{\norm}[2][]{\ifx|#1|\left|\else\csname#1\endcsname|\fi#2\ifx|#1|\right|\else\csname#1\endcsname|\fi}
     \providecommand{\normon}[2]{\norm{#1}_{#2}}

     \providecommand{\abs}[2][]{\ensuremath{\ifx|#1|{\left|#2\right|}\else{\csname#1\endcsname|{#2}\csname#1\endcsname|}\fi}}

     \providecommand{\Norm}[2][]{\ifx|#1|\left\|\else\csname#1\endcsname\|\fi{#2}\ifx|#1|\right\|\else\csname#1\endcsname\|\fi}

     \providecommand{\Normon}[2]{\Norm{#1}_{#2}}
     \providecommand{\Normonspace}[2]{\Norm{#1}_{\vecspace{#2}}}
     
     \providecommand{\normonsob}[4][]{\normon{#2}{\sob{#3}{#4}\if|#1|{}\else(#1)\fi}}
     \providecommand{\Normonsob}[4][]{\Normon{#2}{\sob{#3}{#4}\if|#1|{}\else(#1)\fi}}
     \providecommand{\Normonleb}[4][]{\Normon{#2}{\leb{#3}\if|#1|{}\else(#1)\fi}}
     \providecommand{\ltwop}[3][]{\ensuremath{\qa{#2,#3}\ifx|#1|\else_{#1}\fi}}
     \providecommand{\ltwopreg}[2]{\ensuremath{\qareg{#1,#2}\ifx|#1|\else_{#1}\fi}}
     \providecommand{\ltwopbig}[2]{\ensuremath{\qabig{#1,#2}\ifx|#1|\else_{#1}\fi}}
     \providecommand{\ltwopBig}[2]{\ensuremath{\qaBig{#1,#2}\ifx|#1|\else_{#1}\fi}}
     \providecommand{\ltwopbigg}[2]{\ensuremath{\qabigg{#1,#2}\ifx|#1|\else_{#1}\fi}}
     \providecommand{\ltwopBigg}[2]{\ensuremath{\qaBigg{#1,#2}\ifx|#1|\else_{#1}\fi}}

     \providecommand{\average}[2][]{{\qa{#2}\ifx|#1|\else_{#1}\fi}}

     \providecommand{\ensemble}[2]{\ensuremath{\left\{ #1:\;#2 \right\}}}
     \providecommand{\setofsuch}{\ensemble}%
     \providecommand{\setof}[1]{{\qc{#1}}}
     \providecommand{\pair}[2]{\qp{#1,#2}}

     \providecommand{\conditionalto}[1]{{\left|{#1}\right.}}

    \providecommand{\measure}[1]{\ensuremath{\mathcalbf{\MakeUppercase{#1}}}}
    
    \providecommand{\probmeasure}[2][]{{\measure{#2}}\ifx|#1|\else_{#1}\fi}
    \providecommand{\Prob}{}
    \renewcommand{\Prob}[1][]{\probmeasure[{#1}]{p}}

    \providecommand{\randvars}[1][\Prob]{\operatorname{RV}\ifx|#1|{}\else{(#1)}\fi}
    \providecommand{\discrandvars}[1][\Prob]{\operatorname{DRV}\ifx|#1|{}\else{({#1)}\fi}} 
    \providecommand{\contrandvars}[1][\Prob]{\ensuremath{\operatorname{CDRV}\ifx|#1|{}\else(#1)\fi}} 
     \def\env@matrix{\hskip -\arraycolsep
      \let\@ifnextchar\new@ifnextchar
      \array{*\c@MaxMatrixCols c}}
     \makeatletter
     \renewcommand*\env@matrix[1][c]{\hskip -\arraycolsep
       \let\@ifnextchar\new@ifnextchar
       \array{*\c@MaxMatrixCols #1}}
     \makeatother
     \providecommand{\irow}[2]{#1_{#2}}%
     \providecommand{\icol}[2]{#1^{#2}}%

     \providecommand{\ijrowcol}[3]{\icol{\irow{#1}{#2}}{#3}}
     \providecommand{\entry}[1]{\qb{#1}}
     \providecommand{\vecentry}[2]{\irow{#1}{#2}}

     \providecommand{\rowof}[1]{\qb{#1}}

     \providecommand{\getentryi}[2]{\irow{\entry{#1}}{#2}}
     
     \providecommand{\getvecentry}[2]{\getentryi{\vec #1}{#2}}

     \providecommand{\dismatof}[2][r]{\begin{bmatrix}[#1]#2\end{bmatrix}}

     \providecommand{\matentry}[3]{\ijrowcol{#1}{#2}{#3}}

     \providecommand{\block}[5]{\ijrowcol{#1}{\ifx#2#3{\rowof{#2}}\else\rowof{{#2}\dotsc{#3}}\fi}{\ifx#4#5{\rowof{#4}}\else\rowof{{#4}\dotsc{#5}}\fi}}
     \providecommand{\colblock}[3]{\getvecentry{#1}{\ifx#2#3{#2}\else\fromto{#2}{#3}\fi}}

     \providecommand{\dismatskeldots}[4]{
       \dismatof[c]{
         #1&\dotsc&#3
         \\
         \vdots & \ddots &\vdots
         \\
         #2&\dotsc&#4
       }
     }
     \providecommand{\dismatcommfromtofromto}[5]{
       \dismatskeldots{#1#2#4}{#1#3#4}{#1#2#5}{#1#3#5}
     }
     \providecommand{\dismatcustfromtofromto}[6][matentry]{
       \dismatcommfromtofromto{\csname#1\endcsname{#2}}#3#4#5#6
     }
     \providecommand{\dismatcustfromtofromto}[6][matentry]{
       \dismatskeldots{%
         \csname#1\endcsname{#2}{#3}{#4}%
       }{%
         \csname#1\endcsname{#2}{#3}{#6}%
       }{%
         \csname#1\endcsname{#2}{#5}{#4}%
       }{%
         \csname#1\endcsname{#2}{#5}{#6}%
       }%
     }%
     \providecommand{\dismatcustfromtofromto}[6][matentry]{
       \dismatof{
         \csname#1\endcsname{#2}{#3}{#4}&\dotsc&\csname#1\endcsname{#2}{#3}{#6}
         \\
         \vdots & \ddots &\vdots
         \\
         \csname#1\endcsname{#2}{#5}{#4}&\dotsc&\csname#1\endcsname{#2}{#5}{#6}
       }
     }

     \providecommand{\dissysaxbdotsnm}[5]{\begin{matrix}[r]%
         \matentry{#1}11\vecentry{#2}1&+\dotsb&+\matentry{#1}1{#5}\vecentry{#2}{#5}
         &
         =
         \ifx|#3|0\else{\vecentry {#3}1}\fi
         \\
         \dotsb
         \\
         \matentry{#1}{#4}1\vecentry{#2}1&+\dotsb&+\matentry{#1}{#4}{#5}\vecentry{#2}{#5}
         &
         =
         \ifx|#3|0\else{\vecentry {#3}{#4}}\fi
     \end{matrix}}
     \providecommand{\seqof}[1]{\qp{#1}}%
     \providecommand{\seqs}[2]{\seqof{#1}_{#2}}
     \providecommand{\sets}[2]{\setof{#1}_{#2}}%
     \providecommand{\seqi}[3][]{\seqs{#2_{#3}}{\ifx|#1|{#3}\else{{#3}\in{#1}}\fi}}%

     \providecommand{\seti}[3][]{\sets{#2_{#3}}{\ifx|#1|_{#3}\else_{{#3}\in{#1}}\fi}}%
     \providecommand{\sequ}[3][]{\seqs{#2^{#3}}{\ifx|#1|{#3}\else{{#3}\in{#1}}\fi}}%
     \providecommand{\setu}[3][]{\sets{#2^{#3}}{\ifx|#1|{#3}\else{{#3}\in{#1}}\fi}}%
     \providecommand{\limofat}[3][]{\ensuremath{\lim_{\ifx|#1|{}\else{#1\ni}\fi#3}{#2}}}
     \providecommand{\limsupofat}[3][]{\ensuremath{\limsup_{\ifx|#1|{}\else{#1\ni}\fi#3}{#2}}}
     \providecommand{\liminfofat}[3][]{\ensuremath{\liminf_{\ifx|#1|{}\else{#1\ni}\fi#3}{#2}}}

     \providecommand{\jump}[2][]{\ensuremath{\left\llbracket #2\right\rrbracket\ifx|#1|{}\else_{#1}\fi}}

     \providecommand{\fromto}[2]{\ensuremath{\setof{#1\dotsc#2}}}%

     \providecommand{\d}{}
     \renewcommand{\d}[1][]{\ensuremath{\operatorname{d}\!\ifx|#1|\else{_{#1}}\fi}}
     
     \providecommand{\ds}[1][]{\d{\measure S}}
     \providecommand{\D}[1][]{\ensuremath{\operatorname{D}\!\ifx|#1|\else{_{#1}}\fi}}

    \providecommand{\registered}%
    {\ensuremath{^\text{\textregistered}}}

    \providecommand{\tand}{\ensuremath{\text{ and }}}

    \providecommand{\tor}{\ensuremath{\text{ or }}}

    \providecommand{\constant}[1]{\ensuremath{C_{#1}}}
    \providecommand{\constext}[2][]{\constant{\textup{#2}{\ifx|#1|{}\else{,\ensuremath{#1}}\fi}}}            %
    \providecommand{\constref}[2][]{\ensuremath{\constant{\textup{\ref{#2}{\ifx|#1|{}\else{,\ensuremath{#1}}\fi}}}}}
    \providecommand{\constdef}[2][]{\label{#2}\ensuremath{\constant{\textup{\ref{#2}{\ifx|#1|{}\else{,\ensuremath{#1}}\fi}}}}}

    \providecommand{\funkref}[3][]{\ensuremath{{#3}_{\textup{\ref{#2}{\ifx|#1|{}\else{,\ensuremath{#1}}\fi}}}}}

    \providecommand{\Cof}{\operatorname{Cof}}

    \providecommand{\diam}{\operatorname{diam}}
    
    \providecommand{\curl}{\operatorname{curl}}
    \renewcommand{\div}[1][]{\nabla\ifx|#1|{}\else\kern-2pt_{#1}\fi\kern-2pt\inner}
    \providecommand{\divof}[2][]{\div[#1]\ifx|#2|{}\else\qb{#2}\fi}
    \providecommand{\grad}{}
    \renewcommand{\grad}[1][]{\nabla\ifx|#1|\else_{#1}\fi}

    \providecommand{\rot}[1][]{\nabla\ifx|#1|\else_{#1}\fi\outerp}
    
    \providecommand{\rowdiv}[1][]{\D\ifx|#1|{}\else\kern-1pt_{#1}\kern-2pt\fi\cdot}
    \providecommand{\rowdivof}[2][]{\rowdiv[#1]\ifx|#2|{}\else\qb{#2}\fi}

    \providecommand{\area}{\measure s}

    \providecommand{\inverse}[2][]{\powp[#1]{-1}{#2}}

    \providecommand{\fracl}[3][]{\ifx|#1|\nicefrac{#2}{#3}\else{#2}#1/{#3}\fi}

    \providecommand{\qpfracl}[3][]{\qp{\ifx|#1|\fracl{#2}{#3}\else{#2}#1/{#3}\fi}}
    \providecommand{\qpfrac}[3][]{\qp{\ifx|#1|\frac{#2}{#3}\else{#2}#1/{#3}\fi}}
    \providecommand{\absfracl}[3][]{\abs{\ifx|#1|\fracl{#2}{#3}\else{#2}#1/{#3}\fi}}
    \providecommand{\absfrac}[3][]{\abs{\ifx|#1|\frac{#2}{#3}\else{#2}#1/{#3}\fi}}
    \providecommand{\fraclff}[3][]{\ifx|#1|{#2}/{#3}\else{#2}#1/{#3}\fi}

    \providecommand{\eye}[1][]{\vec{\mathrm I}\ifx|#1|{}\else_{#1}\fi}%
    \providecommand{\numeye}[1][]{\boldsymbol{\mathsf{I}}\ifx|#1|{}\else_{#1}\fi}%
    \providecommand{\Eye}[1]{
      \begin{bmatrix}
      \ifthenelse{#1>1}{
        \ifthenelse{#1>2}{
          \ifthenelse{#1>3}{
            \ifthenelse{#1>4}{
              1&\zeroentry&\dotso&\zeroentry
              \\
              \zeroentry&1&\dotso&\zeroentry
              \\
              \vdots&\vdots&\ddots&\vdots
              \\
              \zeroentry&\zeroentry&\dotso&1
            }{        
              1&\zeroentry&\zeroentry&\zeroentry
              \\
              \zeroentry&1&\zeroentry&\zeroentry
              \\
              \zeroentry&\zeroentry&1&\zeroentry
              \\
              \zeroentry&\zeroentry&\zeroentry&1
            }
          }{
            1&\zeroentry&\zeroentry
            \\
            \zeroentry&1&\zeroentry
            \\
            \zeroentry&\zeroentry&1
          }
        }{
          1&\zeroentry
          \\
          \zeroentry&1
        }
      }{
        1
      }
      \end{bmatrix}
    }

    \providecommand{\lebmeas}[1][]{\operatorname{l}^{#1}}     %
    \providecommand{\lebmeasof}[2][]{\ifx|#1|\left|#2\right|\else\lebmeas[#1]\qp{#2}\fi}         %
    \providecommand{\area}[1]{\operatorname{area}#1}          %

    \providecommand{\meshsize}[1][]{h\ifx|#1|\else_{#1}\fi}

    \providecommand{\mini}[2]{#1\wedge#2}                     %

    \providecommand{\argmin}{\operatorname{argmin}}

    \providecommand{\Argmax}{\operatorname{Argmax}}

    \providecommand{\dash}[1][']{\ifthenelse{\equal{#1}{'}\OR\equal{#1}{''}}{#1}{^{(#1)}}}
    
    \providecommand{\pdfrac}[2][]{\ensuremath{\frac{\partial\ifx|#1|\phantom{#2}\else{#1}\fi}{\partial{#2}}}} %
    \providecommand{\pdfracpow}[3][]{\ensuremath{\frac{\partial^{#3}\ifx|#1|\phantom{#2}\else{#1}\fi}{\partial{#2}^{#3}}}} %

    \providecommand{\pd}[2][]{\ensuremath{\partial_{#2}}{\ifx|#1|{}\else{\qb{#1}}\fi}} %

    \renewcommand{\Im}{\operatorname{im}}                 %
    \renewcommand{\Re}{\operatorname{re}}                 %
    \providecommand{\imaginpart}[1][]{\Im{\ifx|#1|{}\else\qp{#1}\fi}} %
    \providecommand{\realpart}[1][]{\Re{\ifx|#1|{}\else\qp{#1}\fi}} %
    \providecommand{\trace}{\operatorname{tra}}             %

    \providecommand{\transpose}{\intercal}%

    \providecommand{\transposed}{{}^\transpose}

    \providecommand{\orthogonalto}[1][]{\ensuremath{\perp\ifx|#1|{}\else{_{#1}}\fi}}
    
    \providecommand{\rowof}[1]{\ensuremath{\vecof{#1}}}

    \providecommand{\zeroentry}{\phantom0}%

    \providecommand{\smint}{\ensuremath{{\text{\textbf{/}}}\kern-.75em\smallint}}
    \renewcommand{\smint}[1][]{\lower12.3pt\hbox{\begin{tikzpicture}\draw[line width=.75pt] (-3pt,-0.5)--(1pt,-0.5) node[pos=0.6]{$\int$};\path (-1.5pt,-24pt)node {\scriptsize $#1$};\end{tikzpicture}}}

    \providecommand{\lap}{\ensuremath{\Delta}}
    \providecommand{\lapin}[1][]{\lap\ifx|#1|\else_{#1}\fi}
    \providecommand{\normalsymbol}{\operatorname{n}}
    \providecommand{\normal}[1][]{\normalsymbol\ifx|#1|\else_{#1}\fi}%
    \providecommand{\normalto}[1]{\ensuremath{\normal[#1]}}
    \providecommand{\normalder}[1][]{\ensuremath{\normal\ifx|#1|\else\qp{#1}\fi{\inner\grad}}}
    \providecommand{\normalderto}[2][]{\ensuremath{\normalto{#2}\ifx|#1|\else\qp{#1}\fi{\inner\grad}}}

    \providecommand{\tangentialsymbol}{\operatorname{t}}
    \providecommand{\tangentialto}[2][]{\tangentialsymbol\ifx|#1|\else^{#1}\fi\ifx|#2|\else_{#2}\fi}

    \providecommand{\intersected}{\ensuremath{\cap}}
    
    \providecommand{\meet}{\intersected}

    \ifthenelse{\boolean{usesvjour}}{}{
      \let\vec\undefined
      \providecommand{\vec}[1]{\ensuremath{\boldsymbol{#1}}}
      \renewcommand{\vec}[1]{\ensuremath{\boldsymbol{#1}}}
    }

    \providecommand{\hatmat}[1]{\hat{\mat{#1}}}

    \providecommand{\geomat}[1]{\vec{\MakeUppercase{#1}}}

    \providecommand{\mat}[1]{\geomat{#1}} %
    \providecommand{\Prob}[1][]{\ensuremath{\operatorname{Prob}\ifx|#1|{}\else_{#1}\fi}}

    \providecommand{\pdf}[2][]{\ensuremath{\operatorname{pdf}_{#2\ifx|#1|{}\else{\conditionalto{#1}}\fi}}\xspace}

    \providecommand{\expectation}{\ensuremath{\operatorname{E}}}
    \providecommand{\EX}[1][]{\ensuremath{\expectation\ifx|#1|{}\else_{#1}\fi}}

    \providecommand{\gausskernel}[3][x]{%
      \ensuremath{
        \exp\frac{-\if#20{#1}\else(#1-\mu)\fi^2}{%
          2\if#31{}\else\powp2{#3}\fi}%
      }%
    }
    \providecommand{\gaussdistribution}[3][x]{%
      \ensuremath{\frac1{\sqrt{2\pic}\if#31{}\else#3\fi}%
        \gausskernel[#1]{#2}{#3}
      }%
    }%

    \providecommand{\boundary}{\partial}

    \providecommand{\SPD}{\operatorname{SPD}}
    \providecommand{\spdmats}[2][F]{\SPD(\csname#1\endcsname{#2})}
     \providecommand{\Continuous}{\ensuremath{\operatorname C}\xspace}%
     \providecommand{\Hspace}{\ensuremath{\operatorname H}\xspace}
     \providecommand{\Lebesgue}{\ensuremath{\operatorname L}\xspace}

     \providecommand{\Weaklyder}{\ensuremath{\operatorname W}\xspace}
     
     \providecommand{\dual}[1]{\ensuremath{{#1}'}}
     
     \providecommand{\dualspace}[2][]{\dual{\linspace{#2}\ifx|#1|\else{_{#1}}\fi}}
     \providecommand{\bidual}[1]{\ensuremath{{#1}''}}
     \providecommand{\bidualspace}[2][]{\bidual{\linspace{#2}\ifx|#1|\else{_{#1}}\fi}}

     \providecommand{\cont}[1]{\ensuremath{\Continuous^{#1}}}

     \providecommand{\BV}[1]{\ensuremath{\operatorname{BV}}}

     \providecommand{\leb}[1]{\ensuremath{\Lebesgue_{#1}}}
     \providecommand{\lebloc}[1]{\ensuremath{{{\Lebesgue}^{\kern-.20em\lower .1ex\hbox{\tiny\textrm{\textup{loc}}}}_{#1}}}}
     \providecommand{\lebnorm}[3][]{\ensuremath{\Norm{#2}_{\leb{#3}\ifx|#1|{}\else(#1)\fi}}}

     \providecommand{\sob}[2]{\ensuremath{{\smash\Weaklyder}^{#1}_{#2}}}

     \providecommand{\sobh}[1]{\ensuremath{\Hspace^{#1}}}
     \providecommand{\vecsobh}[1]{\ensuremath{\vec\Hspace^{#1}}}
     \providecommand{\hdiv}[1][]{\vecsobh{\operatorname{div}}\ifx|#1|\else(#1)\fi}
     \providecommand{\hcurl}[1][]{\vecsobh{\operatorname{curl}}\ifx|#1|\else(#1)\fi}
     
     \providecommand{\sobhz}[2][]{\sobh{#2}_{0\ifx+#1+\else|#1\fi}}

     \providecommand{\Lip}[1][]{\ensuremath{\operatorname{Lip}}\ifx|#1|{}\else{\qp{#1}}\fi}

     \providecommand{\poly}[1]{\ensuremath{\rP}^{#1}}

     \providecommand{\Sym}{\operatorname{Sym}}
     \providecommand{\Symmatrices}[2][R]{\ensuremath{\operatorname{Sym}{(\csname#1\endcsname{#2})}}}
     
     \providecommand{\SAmatrices}[2][F]{\ensuremath{\operatorname{SA}{(\csname#1\endcsname{#2})}}}
     \providecommand{\mesh}[2][]{{\ensuremath{\mathcalbf{\MakeUppercase{#2}}\ifx|#1|\else_{#1}\fi}}}

    \providecommand{\crouzeixraviart}[1][1]{\operatorname{CR}\ifx|#1|{}\else{^{#1}}\fi}

    \providecommand{\linspace}[1]{\mathscript{\MakeUppercase{#1}}}
    
    \providecommand{\vecspace}{\linspace}
    \providecommand{\linop}[1]{\mathcalbf{\MakeUppercase{#1}}}
    \providecommand{\nlop}[1]{\mathcal{\MakeUppercase{#1}}}

    \providecommand{\Lin}{\operatorname{Lin}}

    \providecommand{\linops}[2]{\ensuremath{\Lin\qp{{#1}\to{#2}}}}
    \providecommand{\linopss}[2]{\linops{\linspace{#1}}{\linspace{#2}}}
    
    \providecommand{\clinopss}[2]{\clinopss{\linspace{#1}}{\linspace{#2}}}
    \providecommand{\fepartition}[2][]{\mathscript{\MakeUppercase{#2}}\ifx|#1|{}\else_{#1}\fi}
    \providecommand{\fespace}[2][]{\mathbb{\uppercase{#2}}\ifx|#1|{}\else_{#1}\fi}
    
    \providecommand{\vespace}[1][]{\fespace v\ifx|#1|\else_{#1}\fi}

    \providecommand{\fe}[2][]{\ensuremath{\uppercase{#2}\ifx|#1|\else_{#1}\fi}}%

    \providecommand{\fespacefun}[2]{\fe[\fespace{#1}]{#2}}
    \providecommand{\vecfe}[2][]{\ensuremath{\vec{\uppercase{#2}}\ifx|#1|{}\else{_{#1}}\fi}}%
    \providecommand{\vecfespacefun}[2]{\vecfe[\fespace{#1}]{#2}}
    \providecommand{\matfe}[2][]{\ensuremath{\mat{\uppercase{#2}}\ifx|#1|{}\else{_{#1}}\fi}}%
    
    \providecommand{\hatmatfe}[2][]{\ensuremath{\hatmat{\uppercase{#2}}\ifx|#1|{}\else{_{#1}}\fi}}%

    \providecommand{\Foreach}{\text{ for each }}%
    \providecommand{\Forsome}{\text{ for some }}

    \RequirePackage{amscd}

    \providecommand{\ideq}{\equiv}
    
    \providecommand{\funk}[3]{\ensuremath{#1:#2\to#3}}
    \providecommand{\ffunk}[3]{\ensuremath{#1:#2\rightrightarrows#3}}

    \providecommand{\implies}{\ensuremath{\:\Rightarrow\:}\xspace}
    \renewcommand{\implies}{\ensuremath{\:\Rightarrow\:}\xspace}

    \providecommand{\restriction}[2]{\left.#1\right|_{#2}}
    \renewcommand{\restriction}[2]{\left.#1\right|_{#2}}

    \providecommand{\evalat}[3][]{\qb{#2}_{\ifx|#1|{}\else#1=\fi#3}}
    \providecommand{\evaldiff}[4][]{\qb{#2}^{\ifx|#1|{}\else#1=\fi#3}_{\ifx|#1|{}\else#1=\fi#4}}
    \providecommand{\boundarytraceof}[2]{\restriction{#2}{\boundary{#1}}}

    \providecommand{\codename}[1]{\nolinkurl{#1}\xspace}

    \ifthenelse{\boolean{useutopia}}{%
      
    }{%
      
    }
    \providecommand{\indexen}[2][]{{\ifthenelse{\boolean{shownotes}}{\color b}{}#2\ifx|#1|\index{#2}\else\index{#1}\fi}}
    \providecommand{\indexemph}[2][]{\emph{\indexen[#1]{#2}}}

    \providecommand{\ListParameters}{}
    \renewcommand{\ListParameters}%
    {
    	 \setlength{\topsep}{0pt}
    	 \setlength{\leftmargin}{0pt}
             \setlength{\itemsep}{0pt}
    	 \setlength{\parsep}{0pt}
    	 \setlength{\parskip}{0pt}
             \setlength{\labelsep}{0pt}
    	 \setlength{\itemindent}{0pt}
    }
    {%
      \begin{list}%
        {}%
        {\ListParameters%
        
    }}%
    {\end{list}}
    \newcounter{tmpcounter}
    \newcounter{LetterListItem}
    \renewcommand{\theLetterListItem}{(\alph{LetterListItem})}
    {
    	\begin{list}%
    	{\theLetterListItem\ }%
    	{\usecounter{LetterListItem}
    	  \ListParameters
              \ifx|#1|{}\else\setcounter{LetterListItem}{#1}\fi
    	}
    }%
    {\end{list}}
    \newcounter{NumberListItem}
    \renewcommand{\theNumberListItem}{\arabic{NumberListItem}}
    \newenvironment{NumberList}%
    {
    	\begin{list}%
    	{\theNumberListItem.\ }%
    	{\usecounter{NumberListItem}%
    	 \ListParameters
    	}
    }%
    {\end{list}}
    \newcounter{QuestionListItem}
    \renewcommand{\theQuestionListItem}{\textbf{Question \arabic{QuestionListItem}}}
    {
    	\begin{list}%
    	{\theQuestionListItem.\ }%
    	{\usecounter{QuestionListItem}%
    	 \ListParameters
    	}
    }%
    {\end{list}}
    \newcounter{RomanListItem}
    \renewcommand{\theRomanListItem}{(\roman{RomanListItem})}
    {
    	\begin{list}%
    	{\theRomanListItem\ }%
    	{\usecounter{RomanListItem}
    	 \ListParameters
    	}
    }%
    {\end{list}}
    \newcounter{StepsItem}
    {
    	\begin{list}%
    	{Step \theStepsItem.\ }%
    	{\usecounter{StepsItem}%
    	 \ListParameters
    	}
    }%
    {\end{list}}
    \newcounter{CasesListItem}
    \renewcommand{\theCasesListItem}{\Alph{CasesListItem}}
    {
    	\begin{list}%
    	{\emph{Case \theCasesListItem.}\ }%
    	{\usecounter{CasesListItem}%
    	 \ListParameters
    	}
    }%
    {\end{list}}
    \newcounter{QAListItem}
    \renewcommand{\theQAListItem}{Q\arabic{QAListItem}:}
    {
    	\begin{list}%
    	{\theQAListItem}%
    	{\usecounter{QAListItem}
    	 \ListParameters
    	}
    }%
    {\end{list}}

    \ifthenelse{\boolean{isthesis}}{%
      \setcounter{secnumdepth}{1}%
    }{
      \setcounter{secnumdepth}{2} %
    }
    \providecommand{\ListParameters}{}
    \renewcommand{\ListParameters}
    {
    	 \setlength{\topsep}{0em}
    	 \setlength{\leftmargin}{0em}
             \setlength{\itemsep}{0ex}
    	 \setlength{\parsep}{.5ex}
    	 \setlength{\itemindent}{\labelsep}
    	 \addtolength{\itemindent}{\labelwidth}
    }

      \providecommand{\ObsName}{Remark}%
      \providecommand{\RemName}{Remark}%
      \providecommand{\NotName}{Notation}%
      \providecommand{\BFNName}{Big~Fat~Note}%
      \providecommand{\DefName}{Definition}%
      \providecommand{\ExaName}{Example}%
      \providecommand{\TheName}{Theorem}%
      \providecommand{\LemName}{Lemma}%
      \providecommand{\ProName}{Proposition}%
      \providecommand{\CorName}{Corollary}%
      \providecommand{\PbmName}{Problem}%
      \providecommand{\HypName}{Hypothesis}%
      \providecommand{\AlgName}{Algorithm}%
      \providecommand{\ExeName}{Exercise}%
      \providecommand{\SolName}{Solution}%
      \providecommand{\ClaName}{Claim}%
      \providecommand{\EsyName}{Essay}%
      \providecommand{\Proofname}{Proof}%
      \providecommand{\Derivename}{Derivation}%
    
    \ifthenelse{\boolean{isthesis}}{%
      \providecommand{\Thecounter}{The}
    }{%
      \providecommand{\Thecounter}{subsection}
    }
    \makeatletter
    \newcommand{\oltikzgetxy}[3]{%
      \tikz@scan@one@point\pgfutil@firstofone#1\relax
      \edef#2{\the\pgf@x}%
      \edef#3{\the\pgf@y}%
    }
    \makeatother
    \providecommand{\pdfformat}[1]{
       \provideboolean{pdfoutput}
       \setboolean{pdfoutput}{#1}%
      \ifthenelse{\boolean{pdfoutput}}{
        \typeout{using pdf}
\usepackage{pdfsync}
        \providecommand{\graphext}{pdf}
        \renewcommand{\graphext}{pdf}
        \providecommand{\graphextex}{pdf_t}
        \renewcommand{\graphextex}{pdf_t}
      }{
        \typeout{using eps}
        \RequirePackage[dvips]{graphicx,xcolor}
        \providecommand{\graphext}{eps}
        \renewcommand{\graphext}{eps}
        \providecommand{\graphextex}{eps_t}
        \renewcommand{\graphextex}{eps_t}
      }
      \RequirePackage{epsfig}
      \RequirePackage{tikz}
      \RequirePackage{rotating}
      \RequirePackage{graphicx}
      \RequirePackage{xcolor}
      \provideboolean{darkcolortheme}
      \definecolor{SussexFlint}{rgb}{.00,.19,.21}
      \definecolor{SussexGrey}{rgb}{.51,.58,.49}
      \definecolor{SussexOrange}{rgb}{.94,.29,.00}
      \definecolor{SussexYellow}{rgb}{1.00,.73,.00}
      \definecolor{SussexRed}{rgb}{.94,.01,.49}
      \definecolor{SussexPurple}{rgb}{.48,.06,.44}
      \definecolor{SussexGreen}{rgb}{.00,.58,.46}
      \definecolor{OmarGreen}{rgb}{.00,.68,.36}
      \definecolor{SussexBlue}{rgb}{.00,.58,.65}
      \definecolor{OmarBlue}{rgb}{.00,.38,.65}
      \colorlet{a}{OmarBlue}%
      \colorlet{b}{SussexOrange}
      \colorlet{c}{SussexGreen}
      \colorlet{d}{SussexPurple}%
      \colorlet{e}{SussexRed}
      \colorlet{f}{SussexYellow}
      \colorlet{g}{white}%
      \colorlet{h}{SussexGrey}%
      \colorlet{i}{black}%
      \colorlet{j}{SussexFlint}
      \colorlet{colora}{a}
      \colorlet{colorb}{b}
      \colorlet{colorc}{c}
      \colorlet{colord}{d}
      \colorlet{colore}{e}
      \colorlet{colorf}{f}
      \colorlet{colorg}{g}
      \colorlet{colorh}{h}
      \colorlet{colori}{i}
      \colorlet{colorj}{j}
      \newcommand{\mausDarkColorTheme}{
        \colorlet{a}{SussexYellow!50!yellow}
        \colorlet{b}{SussexBlue}%
        \colorlet{c}{SussexRed!50!red}
        \colorlet{d}{SussexOrange!50!yellow}
        \colorlet{e}{SussexGreen!50!green}
        \colorlet{f}{SussexPurple!50!magenta}
        \colorlet{g}{black}%
        \colorlet{h}{SussexFlint!50!black}
        \colorlet{i}{white}%
        \colorlet{j}{SussexGrey}
      }
      \ifthenelse{\boolean{darkcolortheme}}{\mausDarkColorTheme}{}
    }
    \providecommand{\solution}{\textbf{\SolName.}\xspace}

     \newcounter{phantombox}[enumi]%
     \provideboolean{showphantoms}
     \renewcommand{\thephantombox}{\Alph{phantombox}}%
     \newcommand{\phantombox}[1]{\stepcounter{phantombox}%
       \ensuremath{\boxed{%
           {\ifthenelse{\boolean{showphantoms}}{#1}{\phantom{#1}}}%
           {\texttt{\tiny\ \colorbox{i!50}{\color g\thephantombox}}
           }%
         }%
       }%
     }

     \provideboolean{hidesolution}
     \newcommand{\consolution}[2][]{
       \ifthenelse{\boolean{hidesolution}}{#1\setboolean{showphantoms}{false}}{%
         {\setboolean{showphantoms}{true}\color{i!50}\par \small {\solution}\ #2\par\ \\[5pt]}}
     }
     \provideboolean{showmarks}
     \providecommand{\showmarks}[1]{%
       \ifthenelse{%
         \boolean{showmarks}}{%
         \marginpar{%
           \tiny [$#1$ mark\ifthenelse{\equal{#1}1}{\phantom{s}}s]}%
       }{}}%

     \newcommand{\condibreak}{\ifthenelse{\boolean{hidesolution}}{\newpage}{}}

     \providecommand{\qeyword}[1]{\index{#1}\ifthenelse{\boolean{shownotes}}{\texttt{\color{e}[#1]}}{}}
     \providecommand{\sourcecite}[2][]{\index{#1}\ifthenelse{\boolean{shownotes}}{\texttt{\color{d}[source: \cite[#1]{#2}]}}{}}

     \RequirePackage{hyphenat}
\newtheoremstyle{plain}%
  {}%
  {}%
  {\mdseries\slshape}%
  {\parindent}%
  {\bfseries}%
  {.}%
  {.5em}%
  {}%

\newtheoremstyle{note}%
  {}%
  {}%
  {}%
  {\parindent}%
  {\bfseries}%
  {.}%
  {.5em}%
  {}%

\newtheoremstyle{claim}%
  {}%
  {}%
  {\mdseries\slshape}%
  {}%
  {\bfseries}%
  {}%
  {.5em}%
  {}%

\newtheoremstyle{exercise}%
  {}%
  {}%
  {}%
  {}%
  {\bfseries}%
  {.}%
  {1em}%
  {}%

\newtheoremstyle{break}%
  {}%
  {}%
  {}%
  {}%
  {\bfseries}%
  {.}%
  {\newline}%
  {}%

\swapnumbers{
  \theoremstyle{plain}
  \ifthenelse
      {\boolean{isthesis}}
      {\newtheorem{The}{\TheName}[section]}%
      {
      }%
{
   \theoremstyle{plain}

   \renewcommand{\Thecounter}{subsection}

   \newtheorem*{The*}{\TheName}
   \newtheorem*{Lem*}{\LemName}
   \newtheorem*{Pro*}{\ProName}
   \newtheorem*{Cor*}{\CorName}
   \newtheorem*{Pbm*}{\PbmName}
   \newtheorem*{Hyp*}{\HypName}
   \newtheorem*{Exe*}{\ExeName}
   \newtheorem*{Txx*}{\ExeName} %
   \newtheorem*{Con*}{Conclusion}
   \newtheorem*{Sum*}{Summary}
 }
 {
   \theoremstyle{claim}

 }
 {
   \theoremstyle{note}

   \newtheorem*{Obs*}{\ObsName}

   \newtheorem*{Def*}{\DefName}
   \newtheorem*{Exa*}{\ExaName}
   \newtheorem*{Alg*}{\AlgName}
 }

 {
   \theoremstyle{break}
 }
}

\newenvironment{The}[1][]{%
  \ifx&#1&%
  \subsection{\TheName\xspace}%
  \else%
  \subsection[#1 theorem]{\TheName\ (#1)}%
  \fi%
  \slshape}{%
  \upshape}

\newenvironment{Pbm}[1][]{\subsection{\PbmName\xspace{\ifx&#1&{}\else{ (#1)}\fi}}\slshape}{\upshape}
\newenvironment{Def}[1][]{\subsection{\DefName\xspace{\ifx&#1&{}\else{ of #1}\fi}}}{}
\newenvironment{Obs}[1][]{\subsection{\ObsName\xspace{\ifx&#1&{}\else{ (#1)}\fi}}}{}

\providecommand{\qed}{\vrule height 5pt depth 0pt width 3pt}
\providecommand{\qqed}{{\raggedright{\ \hfill\qed}}}

\newcounter{passo}

{\par\noindent{\bf \Proofname\ #1}\setcounter{passo}{0}}%
{\qqed\par}
{\par\noindent{\bf \Derivename\ #1}\setcounter{passo}{0}}%
{\qqed\par}
\newenvironment{Proof*}[1][{}]%
{\subsection{\Proofname\ #1}\setcounter{passo}{0}}
{\qqed\par}

    \provideboolean{usebibtex}
    \setboolean{usebibtex}{true}
    \usepackage{natbib}
\usepackage{subfig}
\usepackage{graphics}
\renewcommand{\leq}{\leqslant}
\renewcommand{\geq}{\geqslant}
\renewcommand{\rot}{\nabla\!\times\!}%

\renewcommand{\linop}[1]{\mathcal{\MakeUppercase{#1}}}

\providecommand{\abil}[3][a]{#1\qp{#2\,;\,#3}}
\providecommand{\aqua}[2][a]{\abil[#1]{#2}{#2}}
\ifthenelse{\boolean{useutopia}}{
  
}{
  
}
\providecommand{\feop}[2][h]{\mathcal{\MakeUppercase{#2}}\ifx|#1|\else_{#1}\fi}

\ifthenelse{\boolean{useutopia}}{
  \renewcommand{\fe}[2][]{\ensuremath{\mathsfit{#2}\ifx|#1|\else_{#1}\fi}}%
  \renewcommand{\vecfe}[2][]{\ensuremath{\vec{\mathsfit{#2}}\ifx|#1|{}\else{_{#1}}\fi}}%
  \renewcommand{\matfe}[2][]{\ensuremath{\mat{\mathsfit{\MakeUppercase{#2}}}\ifx|#1|{}\else{_{#1}}\fi}}%
}{
  \renewcommand{\fe}[2][]{\ensuremath{\mathsf{#2}\ifx|#1|\else_{#1}\fi}}%
  \renewcommand{\vecfe}[2][]{\ensuremath{\vec{\mathsf{#2}}\ifx|#1|{}\else{_{#1}}\fi}}%
  \renewcommand{\matfe}[2][]{\ensuremath{\mat{\mathsf{\MakeUppercase{#2}}}\ifx|#1|{}\else{_{#1}}\fi}}%
}
\providecommand{\DNewton}{\operatorname{\mathfrak D}}
\providecommand{\FMA}{\operatorname{\mathscript M}}
\providecommand{\FHJB}{\operatorname{\mathscript B}}%
\providecommand{\aposteriori}{a posteriori\xspace}

\providecommand{\apriori}{{a priori}\xspace}

\renewcommand{\aposteriori}{a posteriori\xspace}

\renewcommand{\apriori}{{a priori}\xspace}

\providecommand{\aanswer}[1]{\textit{Answer:}\ \textup{#1}}
\providecommand{\qquestion}[1]{\item\textit{Question:}\ \texttt{#1}\\}
\provideboolean{includeresponses}
\setboolean{includeresponses}{false}
\setboolean{showchanges}{false}%
\setboolean{shownotes}{false}%
\ifthenelse{\boolean{showchanges}}{
  \hypersetup{pagebackref,colorlinks=true,linkcolor=blue,citecolor=blue}
}{}
\renewcommand{\Ae}[1][]{\ensuremath{\ifx|#1|{\ }\else{\:#1\text-}\fi\operatorname{a.e.}}}
\renewcommand{\eye}{\ensuremath{\vec I}}
\renewcommand{\normalsymbol}{\ensuremath{\vec n}}
\renewcommand{\inner}{\transposed}
\renewcommand{\mini}[2]{\min\setof{#1,#2}}

\providecommand{\ourtitle}{}
\renewcommand{\ourtitle}{%
  A least-squares {Galerkin} gradient recovery method for
  fully nonlinear elliptic equations}
\providecommand{\ourshorttitle}{%
 Least-squares~Galerkin~gradient~recovery~for~nonlinear~{PDEs}} 
\ifthenelse{\boolean{usesvjour}}{
  \title{\ourtitle}
  \runningtitle{\ourshorttitle}
}{
  \title[\ourshorttitle]%
        {\ourtitle}
}
\author{%
  Omar Lakkis%
  \and
  Amireh Mousavi}
\ifthenelse{\boolean{usesvjour}}{
  \institute{%
    O. Lakkis
    \at
    University of Sussex, Brighton, England UK,
    \email{lakkis.o.maths@gmail.com}
    \and
    A. Mousavi
    \at
    Isfahan University of Technology, Isfahan, Iran,
    \email{amireh.mousavi@math.iut.ac.ir}
  }
}{
\newenvironment{theorem}[1][]{\begin{The}[#1]}{\end{The}}
\newenvironment{remark}[1][]{\begin{Obs}[#1]}{\end{Obs}}
\newenvironment{definition}[1][]{\begin{Def}[#1]}{\end{Def}}
\newenvironment{problem}[1][]{\begin{Pbm}[#1]}{\end{Pbm}}
}
\begin{document}
\ifthenelse{\boolean{showchanges}}{%
  \section*{Responses to the referees}
  We thank the anonymous referee for constructive criticism of our paper. Following are our reponses
  and actions taken to address their queries.  The numbering of equations and paragraphs in our
  reponses correspond to the updated manuscript.  We hope the links make reviewing easier.
  \begin{NumberList}
    \qquestion{in the title of section 2: galerkin should be Galerkin}
    \aanswer{Fixed title of \S\ref{sec:a-least-square-for-linear-non-divergence-form}.}
    \qquestion{What are the norms used in the Cordes conditions (4) and (5)?}
    \aanswer{Frobenius. We have added the definition after \eqref{eqn:def:special-Cordes-condition}}
    \qquestion{in eq. (6) I would write $n_\W{n_\W^T}$ instead of
      $n_\W{n_\W}$.}  \aanswer{Ok. We had an extra dot-product sign
      signifying transpose, which is correct but we agree that it may be
      confusing; so we changed it to ``$\inner$'' in
      (\ref{eq:def:tangential-trace-tensor}).  To be consistent with
      this change we have replaced all occurrences of ``$\cdot$'' with
      ``$\inner$''.  The Frobenius product ``$\frobinner$'' and
      ``$\trace$'', which were also missing, is now defined in the text
      after (\ref{eq:def:general-nonlinear-HJB-operator}).  }
    \qquestion{in eq. (11): what is a suitable choice for $\theta$?}
    \aanswer{
      There is no universal optimal value for $\theta$? It's use
      is related to the analysis.  Empirically, most cases
      $\theta=\fracl12$ works well, while the ``pure'' $\theta=0$ or $1$
      choices lead to simpler implementations and work well too.  We
      have a comment on this after (\ref{eq:def:theta-linear-operator}).
      In this paper, we do not discuss the optimal value of $\theta$ and
      $\theta = 0, \dfrac{1}{2}, 1$ are just representatives for
      different values ($\theta \in [0,1]$).  As we see in
      Theorem~\ref{the:convergence-rate}, the convergence order of the
      method is independent of the $\theta$'s value.}
    \qquestion{in eq. (20): I don’t understand the exponent
      $k\wedge\beta$.}
    \aanswer{That means the minimum of $k$ and $\beta$; we have changed this notation
      to $\mini k\beta$ in (\ref{eqn:convergence-rate}).}
    \qquestion{in Definition 3: Amèpre should be Ampère.}
    \aanswer{Fixed.}
    \qquestion{in Definition 4: the supescript α is sometimes missing,
      moreover, I don’t understand the condition $\lambda=0$ (just above
      eq. (29)).}
    \aanswer{We have rephrased parts of
      Definition~\ref{def:Hamilton--Jacobi--Bellman-equation}
      in a clearer way, taking into account this comment.}
    \qquestion{The modification of the finite element method in (35)
      should be presented in section 3 and not in section 4, which is
      about numerical results.}  \aanswer{We have moved this to
      Remark~\ref{obs:modified-method-for-zero-tangential-trace}.}
    \qquestion{Figure 1 and Figure 2 are a bit messy, they contain too
      much number that can be barely read. Moreover, what is the added
      value of Figure 2, it is just the same numerical result with a
      different R. I would suggest to skip one of the two figures.}
    \aanswer{\todo}
    \qquestion{A more fundamental issue is the distinction between
      linear and nonlinear problems. The nonlinear version of the PDE is
      not well represented in the paper. In particular, what is the
      difference between the operators F and F in Definition 2 and
      Theorem 4, respectively.}
    \aanswer{\todo}
    \qquestion{What is the relation between the Hamilton-Jacobi-Bellman
      equation and the Monge-Ampère equation?}
    \aanswer{\todo}
    \qquestion{the paper contains some minor typos.}
    \aanswer{\todo{spell-check}}
  \end{NumberList}

  \setcounter{page}0
  \clearpage
}{}%
\maketitle
\begin{abstract}
  We propose a least squares Galerkin based gradient recovery to
  approximate Dirichlet problems for strong solutions of linear
  elliptic problems in nondivergence form and corresponding \apriori
  and \aposteriori error bounds.  This approach is used to tackle
  \changes{fully nonlinear elliptic problems, e.g., Monge--Ampère,
    Hamilton--Jacobi--Bellman, using the smooth (vanilla) and the
  semismooth Newton linearization.}  We discuss numerical results,
  including adaptive methods based on the \aposteriori error
  indicators.
\end{abstract}
\section{Introduction}
Let $\W$ denote a bounded convex domain in $\R d$, $d\in\naturals$ (typically $d=2,3$). 
\changefromto{
A general elliptic operator $\nlop F$ acting on a function $\funk v\W\reals$ can
be written in the following Bellman form \citep{CaffarelliCabre:95:book:Fully,Krylov:18:book:Sobolev}:
\begin{equation}
  \label{eq:def:general-nonlinear-HJB-operator}
  \begin{gathered}
    \nlop F[v](\vec x)=\sup_{\alpha\in\cA}\qp{\linop L^\alpha v(\vec x)-f^\alpha(\vec x)}
    \\
    \tand
    \linop L^\alpha v(\vec x)
    =
    \mat A^\alpha(\vec x)\frobinner\D^2 v(\vec x)
    +\vec b^\alpha(\vec x)\inner\grad v(\vec x)
    -c^\alpha(\vec x)v(\vec x)-f^\alpha(\vec x),
  \end{gathered}
\end{equation}
where for each $\alpha$ in the given parameter set $\cA$
(corresponding to $\nlop F$) $\mat A^\alpha(\vec x)$ is $\vec
x$-uniformly symmetric positive tensor on $\R d$
\changes{($\mat{M}\frobinner\mat{N}:=\trace\mat M\transposed\mat N$
  for $\mat{M}\in\realmats{n}m$ and, if $m=n$, $\trace\mat{M}$ is the
  trace of $\mat{M}\in\realmats{n}n$ which the sum of eigenvalues of
  $\mat M$)}, $\vec b^\alpha(\vec x)\in\R d$ and $c^\alpha(\vec
x),f^\alpha(\vec x)\in\reals$, and $\grad{v},\D^2 v$ denote the
gradient and the Hessian of $v$. The Dirichlet problem of finding a
function $u$ such that
\begin{equation}
  \label{eq:general-nonlinear-elliptic-Dirichlet-bvp}
  \nlop F[u]=0
  \tand
  \boundarytraceof\W u=r.
\end{equation}
}{
Consider the Dirichlet problem of finding a function $\funk u\W\reals$ such that
\begin{equation}
  \label{eq:general-nonlinear-elliptic-Dirichlet-bvp}
  \nlop F[x, u, \grad{u} ,\D^2 u]=0
  \tand
  \boundarytraceof\W u=r.
\end{equation}
Here,  $\grad{u},\D^2 u$ denote the gradient and the Hessian of $u$ 
and $\funk {\nlop F}{\W \times \reals \times \R d \times \R {d \times d}}{\reals}$ 
is assumed to be elliptic and Newton differentiable which is defined by 
Definition~\ref{def:Newton-differentiable-op}.
}

While viscosity solutions are possible, in a natural way, for this
type of equations, we here focus on smoother solutions. Namely, we
look at the numerical approximations of $u$ in $\sobh2(\W)$ satisfying
(\ref{eq:general-nonlinear-elliptic-Dirichlet-bvp}), termed
\indexemph{strong solution}.  We follow a series of papers on the
matter
\citep{SmearsSuli:16:article:Discontinuous,FengJensen:17:article:Convergent,GallistlSuli:19:article:Mixed},
but with a focus on the different and somewhat more flexible numerical
methodology of least squares gradient recovery Galerkin finite
element method to discretize the linear equations in nondivergence
form that ensue from linearizing
(\ref{eq:general-nonlinear-elliptic-Dirichlet-bvp}) using 
\changefromto{
either a (possibly semismooth) Newton method or the 
Hamilton--Jacobi--Bellman
formulation. This entails two concurrent discretizations for the
spatial variable in $\W$ and for the HJB index in $\cA$.}
{semismooth Newton method.} 
We only state the results here, respectively referring for the details 
of \S\ref{sec:a-least-square-for-linear-non-divergence-form} and
\S\ref{sec:linearization-of-fully-nonlinear} to
\citet{LakkisMousavi:19:article:A-least-squares} and
\citet{LakkisMousavi:20:unpublished:A-least-squares}.  We look at some
numerical examples, outlining an adaptive algorithm based on a
posteriori error estimates for the linear elliptic equations in
nondivergence form with Cordes coefficients.
\section{A least-squares \changes{Galerkin} approach to gradient recovery for linear
  equations in nondivergence form}
\label{sec:a-least-square-for-linear-non-divergence-form}
\changes{We outline the proposed numerical method of the strong solution of the
linear second order equation in nondivergence form; the details,
including the proofs of all stated results can be found
in~\cite{LakkisMousavi:19:article:A-least-squares}.} To prevent
difficulties arising from numerically working in $\sobh2(\W)$ space,
we consider an equivalent problem with solution in a
$\sobh1$-regularity space. For this we minimize a \indexemph{cost}
(least-squares) functional associated to the main problem.  We prove
that the equivalent problem is well posed using a coercivity argument,
deducing thus the same result for the discrete counterpart.  By
setting Galerkin finite element spaces within $\sobh1(\W)$, we provide
\apriori and \aposteriori error bounds.

Dropping the index $\alpha$ from the $\linop L^\alpha$ in
(\ref{eq:def:general-nonlinear-HJB-operator}), we consider the
following linear second order elliptic equations in nondivergence form of
finding $u\in\sobh2(\W)$ such that
\begin{equation}
  \label{eq:nondivergence-homogeneous}
  \linop Lu
  :=
  \mat A\frobinner\D^2 u
  +
  \vec b\inner \nabla u
  -
  c u
  =
  f
  \tand
  \restriction u{\boundary\W}=0
\end{equation}
where
the coefficients
$\mat{A}\in\leb\infty(\W;\Symmatrices{d})$, with
$\Sym(X)=:\text{\indexen{symmetric~operators on} X}$,
is uniformly elliptic, $\vec{b}\in\leb\infty(\W;\R d)$ and
$c\in\leb\infty(\W)$, $c\geq0$ satisfy exactly one of the following
two \changes{\indexemph{Cordes conditions} for some $\varepsilon\in\opinter01$
\begin{gather}
  \label{eqn:def:general-Cordes-condition}
  \vec b\neq\vec0\tor c\neq0
  \implies
  \dfrac{\norm{\mat A}^2
    + \fracl{\norm{\vec b}^2}{2\lambda}
    + (\fracl {c}\lambda)^2}{
    (\trace\mat A
    + \fracl {c}\lambda ) ^2
  }
  \leq
  \dfrac{1}{d+\varepsilon}
  \text{\Ae in } \W
  \Forsome\lambda>0,,
  \\
  \label{eqn:def:special-Cordes-condition}
  \vec b \ideq\vec 0\tand c \ideq0
  \implies
  \dfrac{\norm{ \mat A}^2}{(\trace \mat A)^2}
  \leq
  \dfrac{1}{d-1+\varepsilon}
  \text{\Ae in }\W
  ,
\end{gather}
where $\norm{\mat X}=\powsqrt{\big(\trace\mat X\inner\mat{X} \big)}$.
}%
The right-hand side $f$ is a generic element of
$\leb2(\W)$. We consider the right-hand side $r$ in
(\ref{eq:general-nonlinear-elliptic-Dirichlet-bvp}) to be $0$ for
simplicity (although the developments can be extended to
$\boundarytraceof\W r$ being the trace of a function
$r\in\sobh{2}(\W)$.

Problem (\ref{eq:nondivergence-homogeneous}) is well posed under these
assumptions as shown by \cite{SmearsSuli:14:article:Discontinuous}. In
numerical approximating solutions, dealing with more regular than
$\sobh1(\W)$ spaces leads to complicated computations.  To avoid this
difficulty, we intend to consider an alternative equivalent problem
with $\sobh1(\W)$ solution.

\changes{We denote the outer normal to $\W$ at $\vec
x\in\boundary\W$ by $\normalto{\W}(\vec x)$, which we assume defined
for $\area$-almost every $\vec x\in\boundary\W$ ($\area$ being the
$(d-1)$-dimensional ``surface'' measure) and recall the tangential
trace of $\vec \psi \in \sobh1(\W; \R d)$ is expressed (or defined)
by
\begin{equation}
  \label{eq:def:tangential-trace-tensor}
  \qp{\eye - \normalto{\W}\normalto{\W}\inner}\boundarytraceof\W{\vec\psi}. 
\end{equation}}%
Define the following function spaces
\index{$\linspace w$}
\index{$\linspace v$}
\index{$\linspace y$}
\begin{gather}
  \linspace W:= \left\lbrace 
  \vec{\psi}\in\sobh1(\W;\R d) :
  \qp{\eye-\normalto\W\normalto\W\transposed}
  \restriction{\vec\psi}{\boundary\W} = 0 \right\rbrace 
  ,
   \\
  \linspace{Y}:=
  \sobh1(\W) 
  \times
  \sobh1\qp{\W;\R d}
  \\
  \linspace{V}
  :=
  \sobhz1(\W) \times  \linspace W 
  \subseteq
  \linspace Y
  ,
\end{gather}
endowed with the $\sobh1$-norm for $\linspace{W}$ and the following 
norm for $\linspace{Y}$ and $\linspace{V}$,
\begin{equation}
  \Norm{(\varphi , \vec\psi) }_{\linspace{Y}}^2
  :=
  \Norm{ \varphi  }_{\sobh1(\W)}^2
  +
  \Norm{ \vec\psi  }_{\sobh1(\W)}^2
  \Foreach
  (\varphi ,\vec\psi)
  \in
  \linspace{Y}
  \supseteq
  \linspace V
 .
\end{equation}
We denote by $\ltwop\varphi\psi$ the $\leb2(D;V)$ inner product with
respect to the Lebesgue or surface measure on $D$.
For a fixed $\theta\in\clinter01$ we introduce the linear operator
$\funk{\linop{M}_\theta}{\linspace Y}{\leb2(\W)}$
\begin{equation}
  \label{eq:def:theta-linear-operator}
  \begin{gathered}
    (\varphi, \vec \psi)
    \mapsto
    \mat A \frobinner \D \vec \psi
    + 
    \vec b
    \inner
    (
    \theta\vec \psi +(1-\theta) \grad \varphi
    )
    -
    c \varphi
    =:\linop{M}_\theta(\varphi,\vec\psi)
    .
  \end{gathered}
\end{equation}
\changes{%
  The parameter $\theta$ is at the user's disposal, but the
  most useful values are $0$, $\fracl12$ and $1$.
  We introduce the following quadratic functional of $(\varphi,\vec\psi)\in\linspace V$
}
\index{$E_\theta$}
\begin{equation}
  \begin{gathered}
    \label{functional:minimize}
    E_\theta(\varphi,\vec\psi)
    :=
    \Norm{ \grad\varphi-\vec\psi}_{\leb2(\W)} ^2
    +
    \Norm{ \rot\vec\psi }_{\leb2(\W)}^2
    +
    \Norm{ \linop{M}_\theta(\varphi, \vec\psi)  -f } _{\leb2(\W)}^2
  \end{gathered}
\end{equation} 
where $\rot\vec\psi$ denotes curl (rotation) of $\vec \psi$, and then
consider the convex minimization problem of finding
\begin{equation}
  \label{eq:minimization}
  (u, \vec g)
  =\underset
  {\substack{
   \pair\varphi{\vec\psi}\in\linspace V
  }}
  \argmin%
  E_\theta(\varphi ,\vec\psi)
  .
\end{equation}
\begin{remark}[Equivalent problems]
The problem of finding strong\margnote{strong or weak?}
\margnote[amireh]{we are looking for $u \in \sobh2(\W) \meet \sobhz1(\W)$ 
to (\ref{eq:nondivergence-homogeneous}) which is strong solution.}
 solution to
(\ref{eq:nondivergence-homogeneous}) and convex minimization problem
(\ref{eq:minimization}) are equivalent\margnote{do you need to assume
  twice differentiable $u$ for this equivalence to hold?}
  \margnote[amireh]{Without putting such assumption, twice differentiability of $u$ is implicitly concluded. 
  Because we consider $\vec g =\grad u \in \sobh1(\W)$.}
   and in
(\ref{eq:minimization}), $\vec g = \grad u$ holds. Thus, in the rest
of the paper, $\vec g$ is equal to $\grad u$.
\end{remark}

The Euler--Lagrange equation of the minimization problem 
(\ref{eq:minimization}) consists in finding\index{$\vec g$}
\({
  (u, \vec g) \in \linspace V%
}\)
such that
\begin{multline}
  \label{eq:Euler-Lagrange}
  \qa{ \grad u-\vec g, \grad \varphi -  \vec \psi }
  +
  \qa{\rot\vec g,\rot\vec\psi}
  +
  \qa{ \linop{M}_\theta(u, \vec g) , \linop{M}_\theta(\varphi, \vec \psi) }
  \\
  = 
  \qa{ f, \linop{M}_\theta(\varphi, \vec \psi) }
  ~\Foreach (\varphi, \vec\psi) \in \linspace V.
\end{multline}
We introduce the symmetric bilinear form
\({
  \funk{a_\theta}{
    \ppow{
      \linspace{Y} 
    }2
  }
  \reals
}\)
via
\begin{equation}
  a_\theta
      (\varphi ,\vec\psi
    ;
      \varphi' , \vec\psi')
    := 
    \ltwop{
      \grad \varphi-\vec\psi
    }{
      \grad \varphi'-\vec\psi'
    }
    +
    \ltwop{
      \rot\vec\psi
    }{
      \rot\vec\psi'
    }
    +
    \ltwop{
      \linop{M}_\theta(\varphi ,\vec\psi )
    }{
      \linop{M}_\theta(\varphi',\vec\psi')
    }.
\end{equation}
\begin{theorem}[Coercivity and continuity]
\label{the:corecivity-continuity}
Let $\W$ be a bounded convex open subset of $\R d$ and the uniformly 
bounded coefficients $\mat A, \vec b, c$ satisfy either (\ref{eqn:def:general-Cordes-condition}) 
with $\lambda>0$ or (\ref{eqn:def:special-Cordes-condition}) with $\vec b\ideq\vec 0$ 
and $c\ideq0$. Then $a_\theta$ on $\linspace V$ is coercive 
and continuous, there exist $\constref{ineq:coercivity},\constref{ineq:continuity}>0$ 
such that 
 \begin{gather}
 \label{ineq:coercivity}
   \aqua[a_\theta]{\varphi, \vec \psi}
   \geq 
   \constref{ineq:coercivity}
    \Norm{(\varphi, \vec \psi)}_{\linspace y}^2
   \Foreach (\varphi, \vec \psi) \in \linspace V,
   \\
   \begin{aligned}
   \label{ineq:continuity}
     \abil[a_\theta]{\varphi, \vec \psi}{\varphi', \vec \psi'}
   \leq
   \constref{ineq:continuity}
   \Norm{(\varphi, \vec \psi)}_{\linspace y} 
   \Norm{(\varphi', \vec \psi')}_{\linspace y}
   \Foreach (\varphi, \vec \psi), (\varphi', \vec \psi') \in \linspace V.
   \end{aligned} 
\end{gather} 
\end{theorem}
Theorem~\ref{the:corecivity-continuity} ensures the well-posedness of 
the problem~(\ref{eq:Euler-Lagrange}) trough the Lax-Milgram setting.
\begin{definition}[A least squares finite element method]
Let $\mathfrak T$ \index{$\mathfrak T$} be a collection of conforming 
shape-regular triangulations on $\W$ which also known as meshes. If the 
domain, $\W$, is a polyhedral then it coincides with the interior area of 
the mesh. Otherwise, if the domain includes curved boundary, the coincidence 
is lost. Hence this leads to have simplices with curved sides and isoparametric 
elements. For each element $K\in\mathcal{T} \in\mathfrak T$, denote 
$h_{K}:=\diam{K}$, and $h:= h_{\mathcal{T}}:=\max_{{K}\in \mathcal{T}}h_{K}$. 
Now, consider the following Galerkin finite element spaces
\begin{gather}
  \label{def:interpolation-space}
  \fespace u
  :=
  \poly{k}\qp{\mathcal{T}}\meet\sobhz1(\W),
  \quad
  \fespace g
  :=
  \poly{k}\qpreg{\mathcal{T};\R d}\meet \linspace W
  \subseteq
  \sobh1(\W; \R d)
  .
\end{gather}
Corresponding to these spaces, the discrete problem corresponding to
(\ref{eq:Euler-Lagrange}) turns to finding
$(\fespacefun{u}u,\vecfespacefun{g}g)\in\fespace{u}\times\fespace{g}$
such that
\begin{equation}
  \label{eq:discrete}
  a_\theta(\fespacefun uu , \vecfespacefun gg ; \varphi, \vec \psi) 
  =
  \qa{
    f , \linop{M}_\theta(\varphi , \vec \psi)
  }
  \Foreach (\varphi, \vec\psi)\in \fespace u \times \fespace g.
\end{equation}
The coercivity is inherited to subspaces, therefore the solution of discrete 
problem~(\ref{eq:discrete}) is also well-posed.
The discrete problem~(\ref{eq:discrete}) leads to an approximate solution 
satisfying the following error estimate theorems.
\end{definition}
\begin{remark}[implementing the boundary conditions]
  \changes{
    \label{obs:modified-method-for-zero-tangential-trace}
    Since imposing zero-tangential trace condition to the finite
    element spaces is not trivial. In the implementation we used in
    \S\ref{sec:numerical-experiments} we replace in (\ref{eq:discrete}) the space
    $\fespace{g}:=\poly{k}\qp{\mathcal{T};\R{d}}\meet\linspace{W}\subseteq\sobh1(\W;\R{d})$
    with the larger space
    $\tilde{\fespace{g}}:=\poly{k}\qp{\mathcal{T};\R{d}}\meet\sobh1(\W;\R{d})$.
    \margnote[omar]{Should we add a comment on how to make the analysis work for this? The
    $\curl\psi$ penalty should explain that.}
  }
\end{remark}
\begin{theorem}[a priori error estimate]
\label{the:convergence-rate}
   Let $\mathcal{T} \in \mathfrak{T}$ be a mesh on the polyhedral domain 
   $\W\subseteq\R d$. Moreover assume that the strong solution $u$ of 
   (\ref{eq:nondivergence-homogeneous}) satisfies $u \in \sobh{\beta+2}(\W)$,
   for some real $\beta > 0$.  
   Let $(\fespacefun uu,\vecfespacefun gg) \in \fespace u \times \fespace g$  
   be the finite element solution of (\ref{eq:discrete})
   on the mesh $\mathcal{T}$. Then for some 
   $\constref{eqn:convergence-rate}>0$, independent
   of $u$ and $h$ we have
   \begin{equation}
     \label{eqn:convergence-rate}
     \Norm{(u,\grad u) - (\fespacefun u u , \vecfespacefun g g )}_{\linspace Y} 
     \leq
     \constref{eqn:convergence-rate} h^{\mini k\beta}
     \Norm{u}_{\sobh{k+2}(\W)}.
   \end{equation}
\end{theorem}
\begin{remark}[curved domain]
In the case that $\W$ has a curved boundary we use isoparametric
finite element. A piecewise smooth domain
guarantees an optimal rate error bound using isoparametric finite element
similarly to Theorem~\ref{the:convergence-rate} \citep{Ciarlet:2002:book:FEM}.
\end{remark}
\begin{theorem}[error-residual \aposteriori estimates]
\label{the:residual-error-bound}
Let $(\fespacefun u u, \vecfespacefun g g)$ is the unique solution of the 
discrete problem (\ref{eq:discrete}).
 \begin{enumerate}[(i)]
  \item
    The following \aposteriori residual upper bound holds 
    \begin{equation*}
      \label{eqn:aposteriori-error-bound}
      \begin{aligned}
        &\Norm{ (u, \grad u)-(\fespacefun u u,\vecfespacefun g g)}_{\linspace Y}^2
        \leq
		\constref{ineq:coercivity}^{-1}
        \\
        & \qquad
        \qp{
          \Norm{ \grad\fespacefun u u-\vecfespacefun g g}_{\leb2(\W)}^2 
          +
          \Norm{ \rot \vecfespacefun g g}_{\leb2(\W)}^2
          +
          \Norm{ \linop{M}_\theta(\fespacefun u u, \vecfespacefun g g ) -f }_{\leb2(\W)}^2 
        }.
      \end{aligned}
    \end{equation*} 
  \item
    For any open subdomain $\w \subseteq \Omega$ we have 
    \begin{multline}
      \label{eqn:lower-error-bound}
      \Norm{ \grad {\fespacefun u u}- \vecfespacefun g g}_{\leb2(\w)}^2 
      +
      \Norm{ \rot \vecfespacefun g g }_{\leb2(\w)}^2
      +
      \Norm{ \linop{M}_\theta({\fespacefun u u}, {\vecfespacefun g g}) -f }_{\leb2(\w)}^2  
      \\
      \leq
      \constref[\w]{ineq:continuity}
      \qp{ 
        \Norm{u-\fespacefun uu }_{\sobh1{(\w)}}^2
	+
	\Norm{\grad u-\vecfespacefun g g }_{\sobh1{(\w)}}^2 
      },
    \end{multline}
    where $\constref[\w]{ineq:continuity}$ is the continuity constant of 
    $a_\theta$ restricted to $\w \subseteq \Omega$.
 \end{enumerate}
\end{theorem}
\section{Linearization of fully nonlinear problems}
\label{sec:linearization-of-fully-nonlinear}
In this section, we present the \indexen{Newton differentiability} concept
to operators, which can even include non-smooth operators. This concept
is useful to extend the standard Newton linearization to the problems
with non-smooth operator. We state the convergence analysis of a
linearization method which is based on this concept. We then discuss
linearization of two specific fully nonlinear PDEs, namely
Monge--Ampère and Hamilton--Jacobi--Bellman equations that lead to a
sequence of linear equations in nondivergence form.
\begin{definition}[Newton differentiable operator, \citet{ItoKunisch:08:book:Lagrange}]
\label{def:Newton-differentiable-op}
  Let $\linspace X$ and $\linspace Z$ be Banach
  spaces and let $\mathcal U$ be a non-empty open
  subset of $\linspace X$. An operator $\funk{\nlop
    F}{\mathcal{U}\subset\linspace{X}}{\linspace{Z}}$ is called
  \indexemph{Newton differentiable} at $x \in\mathcal U$ if there
  exists a set-valued map with non-empty images
  $\ffunk{\DNewton\nlop{F}}{\mathcal U}{\linopss X Z}$
  (where the double arrow signifies values in the power set of the right-hand
  side) such that
  \begin{equation}
  \label{def:semismooth}
  \lim_{\Normonspace{e}X \rightarrow 0 }
  \dfrac{1}{\Normonspace{e}X}
  \sup_{\linop D\in\DNewton\nlop F[x]}
  \Normonspace{
    \nlop F[x+e] -\nlop F[x]-\linop D e
  }Z 
  = 0
  \Foreach x\in\mathcal U.
  \end{equation}
  The nonlinear operator $\nlop F$
  is called \indexemph{Newton differentiable on} $\mathcal U$ with Newton
  derivative $\DNewton\nlop F$ if $\nlop F$ is Newton differentiable at $x$, for
  every $x \in \mathcal U$.

  The set-valued map $\DNewton\nlop F[x]$ is single-valued at $x$ if
  and only if $\nlop F$ is Fréchet differentiable and $\DNewton\nlop
  F[x]=\setof{\D\nlop F[x]}$.
\end{definition}
\begin{theorem}[Superlinear convergence]
\label{the:superlinear:convergence}
  Suppose that a nonlinear operator $\nlop F$ is Newton differentiable  in an open neighborhood
  $\mathcal U$ of 
  $x^\ast$, solution of $\nlop F[x]=0$. 
  If for any $x \in U$, the all $D \in \DNewton\nlop F[x]$  are non-singular and $\Norm{\inverse D}$ 
  are bounded, then the Newton iteration 
  \begin{equation} 
    \label{eq:recursive}
    x_{n+1} = x_n - D^{-1}_n \nlop F[x_n], \quad D_n \in \DNewton \nlop F[x_n]
  \end{equation}
  converges superlinearly to $x^\ast$ provided that $x_0$ is sufficiently close to $x^\ast$.
\end{theorem}
\begin{definition}[{The Monge--\changes{Ampère} equation}]
  \label{def:monge-ampere}
  Let $\W\subseteq\R2$ be a bounded convex domain. Consider the Monge--Ampère (MA) equation 
  with Dirichlet boundary condition 
  \begin{equation}
    \label{eq:MA}
    \det \D^2 u = f
    \text{ in } \W,
    \restriction u{\boundary\W}=0
    \tand
    u \text{ is strictly convex in } \W,
  \end{equation}
  where $f \in \leb2(\W)$, $f >0$ . 
  Let $\linspace{K}:=\setofsuch{v\in\sobh2(\W)\meet \sobhz1(\W)}{v\text{ is strictly convex}}$
  and define the operators %
  $\funk{\FMA}{\linspace{K}}{ \leb2(\W)}$ by 
  \begin{equation}
    \label{op:MA}
    \FMA [v] := \det\D^2 v -f
  \end{equation}
  and $\funk{\DNewton \FMA}{\linspace{K}}{\linops{\sobh2(\W) \meet \sobhz1(\W)}{\leb2(\W)}}$ 
  by 
  \begin{equation}
    \label{op:D-MA}
    \DNewton \FMA [v] := \Cof \D^2 v : \D^2. 
  \end{equation}
\end{definition}
\begin{theorem}[superlinear convergence of iterative method to MA equation]
  \label{the:superlinear-convergence-M-A}
  The operator $\FMA $ is Fréchet differentiable and thus Newton differentiable.
  Moreover, if the initial guess $u_0 \in \linspace{K}$
  is close to the exact solution 
  $u \in \sobh2(\W) \meet \sobhz1(\W)$ of (\ref{eq:MA}), then the recursive problem
  \begin{equation}
    \label{eq:recursive-MA}
    \Cof\D^2 u_n: \D^2 u_{n+1} = f - \det \D^2 u_n + \Cof\D^2 u_n : \D^2 u_n 
    ~ \text{ in } \W,
    \tand
    \restriction {u_{n+1}}{\boundary\W}=0
  \end{equation}
  converges with superlinear rate to $u$.
\end{theorem}
\begin{definition}[{Hamilton--Jacobi--Bellman equation}]
  \label{def:Hamilton--Jacobi--Bellman-equation}
  Let $\W$ be a bounded convex domain in $\R d$, $d\in\naturals$ (typically $d=2,3$). 
  Consider the Hamilton--Jacobi--Bellman (HJB) equation with Dirichlet boundary condition 
  \begin{equation}
    \label{eq:HJB}
    \sup_{\alpha \in \mathcal{A}}
    \left( 
    \mat A^\alpha\frobinner\D^2 u
    +
    \vec b^\alpha\inner \nabla u
    -
    c ^\alpha u
    - f^\alpha 
    \right) = 0 ~ \text{ in } \W
    \tand
    \restriction u{\boundary\W}=0
  \end{equation}
  where $\mathcal{A}$ is a compact metric space,  
  $\mat{A}\in\leb\infty(\W;\cont0(\mathcal{A};\Symmatrices{d}))$, 
  $\vec{b}\in\leb\infty(\W;\cont0(\mathcal{A};\R{d}))$,
  $c\in\leb\infty(\W;\cont0(\mathcal{A}))$ and 
  $f \in \leb2(\W;\cont0(\mathcal{A}))$. 
  We suppose $\mat A^\alpha(\vec x)$
  is uniformly elliptic in both $\vec x$ and $\alpha$
  and together with $\vec b^\alpha, c^\alpha$ %
  \changes{meets, for some $\epsilon\in\opinter01$
  the Cordes condition (\ref{eqn:def:general-Cordes-condition}), 
  or
  (\ref{eqn:def:special-Cordes-condition}) if $\vec b^\alpha \ideq\vec 0$, $c^\alpha \ideq 0$,
  independent of $\alpha \in \mathcal{A}$.
  }%
  \margnote[Amireh]{we should say this in a way that it be clear the
    constants are independent of $\alpha$}
  \margnote[omar]{i rephrased it, but what happens if for some
    $\alpha$s $b\ideq\vec0, c\ideq 0$ and for some not?}
  For each $\alpha\in\mathcal A$, define the linear operator %
  \begin{equation}
    \label{op:lin-HJB}
    \linop L^\alpha v
    :=
    \mat A^\alpha\frobinner\D^2 v
    +
    \vec b^\alpha\inner \nabla v
    -
    c ^\alpha v,
  \end{equation}
  \changes{the following set of $\mathcal A$-index-valued maps:}
  \begin{equation}
    \mathcal{Q}:= \left\lbrace 
    q:\W \rightarrow {\mathcal{A}}
    \left\vert
    ~ q \text{ is measurable}
    \right\rbrace
    \right.  
    ,
  \end{equation}
  and the set-valued map $\mathcal{N}$, for $v\in\sobh2(\W)\meet\sobhz1(\W)$, such that
  \begin{equation}
    \mathcal{N}[v]:=
    \setofsuch{
      q \in \mathcal{Q}
    }{
      q(\vec x) \in  
      \Argmax_{\alpha\in\mathcal A}\qp{ \qb{\linop L^\alpha v - f^\alpha} \vec x}
      \Aa\vec x \text{ in } \W
    }.
  \end{equation}
  Now, we define the HJB operator %
  by 
  \begin{equation}
    \label{op:HJB}
    \FHJB[v]:%
    =
    \sup_{q \in \mathcal{Q} } 
    \linop L^q v - f^q ,
  \end{equation}
  and the set-valued map 
  \({
    \ffunk{
      \DNewton\FHJB
    }{
      \sobh2(\W)\meet\sobhz1(\W)
    }{
      \linops{\sobh2(\W)\meet\sobhz1(\W)}{\leb2(\W)}
    }
  }\)
  by
  \begin{equation}
    \label{eq:semismooth-HJB}
    \DNewton\FHJB[v]:= \left\lbrace  \linop L^q :=
    (\mat A^q \frobinner\D^2  + \vec b^q \inner \nabla  - c ^q ) 
    \left\vert ~ q \in \mathcal{N}[v] \right\rbrace \right. . 
  \end{equation}
\end{definition}
\begin{theorem}[superlinear convergence of iterative  method to HJB equation]
  \label{the:superlinear-convergence-HJB}
  The operator $\FHJB$ is Newton differentiable with Newton derivative 
  $\DNewton\FHJB$. Moreover, if the initial guess $u_0$ is close 
  to the exact solution $u \in \sobh2(\W) \meet \sobhz1(\W)$ of (\ref{eq:HJB}), 
  the recursive problem
  \begin{equation}
    \label{eq:recursive-HJB}
    \linop L^{q_n} u_{n+1} = f^{q_n} 
    \text{ in }\W,
    \tand
    \restriction {u_{n+1}}{\boundary\W}=0
  \end{equation}
  where $q_n\in\mathcal{N}[u_n]$, converges with superlinear rate to $u$.
\end{theorem}

To follow (\ref{eq:recursive-MA}) and (\ref{eq:recursive-HJB}), we need to approximate 
a linear problem in nondivergence form in each iteration, which we apply the method 
discussed in \S~\ref{sec:a-least-square-for-linear-non-divergence-form}. 
The convergence of the iterative methods (\ref{eq:recursive-MA}) and (\ref{eq:recursive-HJB}) 
implies that the finite element approximation 
$(\fespacefun u u,\vecfespacefun g g) \in \fespace u \times \fespace g$ achieved via the 
recursive problems also satisfies the error bound of Theorem~\ref{the:convergence-rate} 
and \ref{the:residual-error-bound}. 
\begin{remark}
\label{rem:adaptive-refinement}
The \aposteriori residual bound of Theorem~\ref{the:residual-error-bound} 
can be used as an explicit error indicator to determine a locally refined mesh 
in the adaptive scheme.
\end{remark}
\section{Numerical experiments}
\label{sec:numerical-experiments}
We discuss two numerical tests one for each of Monge--Ampère via
Newton and Hamilton--Jacobi--Bellman problems that demonstrate the
robustness of our method to the fully nonlinear problems.  For both
test problems, the domain, $\W$, is taken to be the unit disk in $\R2$
with center at the origin. The criterion to stop the iteration is
either $\Norm{(\fe u_{n+1},\vecfe g_{n+1})-(\fe u_n,\vecfe
  g_n)}_{\linspace{Y}}<10^{-8}$ or maximum $8$ iterations. In
implementation, we take the parameter $\theta$ of (\ref{eq:discrete})
equal to $0.5$. Both implementations were done by using FEniCS
package.

In the first test problem, the known solution is considered smooth and
we see that the numerical results which obtained on the uniform mesh
confirm the convergence analysis of
Theorem~\ref{the:convergence-rate}.
In the second test problem,
we choose the known solution near singular and test the performance of
the adaptive scheme as mentioned in
Remark~\ref{rem:adaptive-refinement}.
Through comparing the convergence rate by the adaptive with uniform refinement, 
we observe the efficiency of the adaptive scheme. 
\begin{problem}[Monge--Ampère test]
Consider problem~(\ref{eq:MA}) and choose $f$ corresponding to the
exact solution
\begin{equation}
  u(\vec x) = -\sqrt{R^2- x_1^2 - x_2^2} +
  \sqrt{R^2-1},\text{ for a fixed }R > 1.
\end{equation}
As suggested by
\citet{LakkisPryer:13:article:A-finite} the first iterate $\fe u_0$
is the discretization of $u_0$ satisfying
\begin{equation}
  \label{eq:initial-guess-HJB}
  \lap u_0 = 2\sqrt{f}
  ~ \text{ in } \W,
  \tand
  \restriction {u_0}{\boundary\W}=0
\end{equation}
and then we track the recursive problem~(\ref{eq:recursive-MA}). 
We show various error norms of linear ($\poly{1}$) and quadratic ($\poly{2}$) 
finite element approximation for two values $R$ in
Figures \ref{fig:MA-R-sqrt2} and \ref{fig:MA-R-2}.
\begin{figure}[tbh]
  \begin{center}
   \subfloat[$\poly{1}$ elements] 
    {
      \includegraphics[width=.475\linewidth]{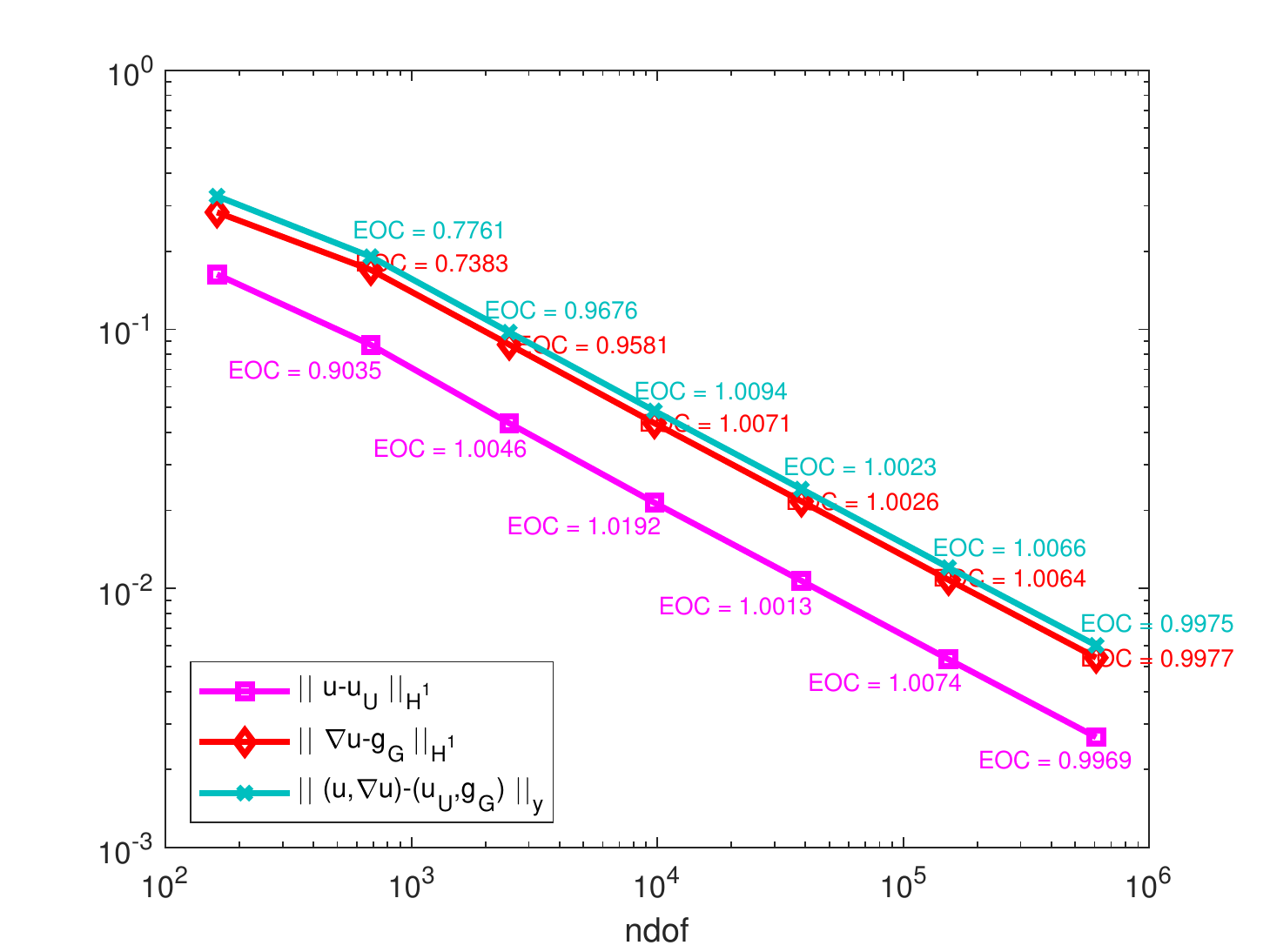}
    }
    \hfill
    \subfloat[$\poly{2}$ elements] 
    {\includegraphics[width=.475\linewidth]{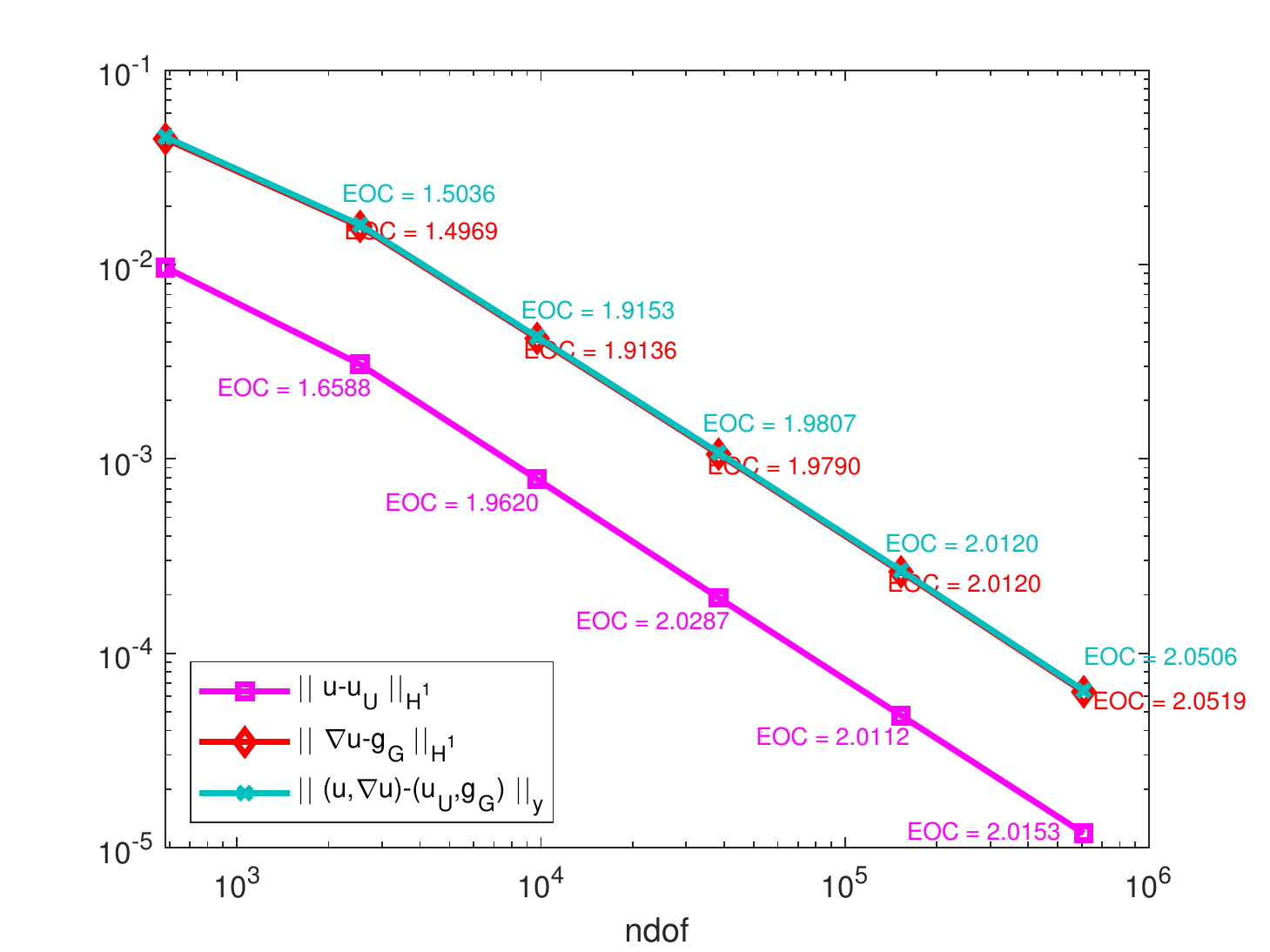}}
    \caption{%
      Experimental order of convergence (EOC) for the Monge--Ampère
      test problem with $R=\sqrt{2}$.}
  \label{fig:MA-R-sqrt2}
  \end{center}
\end{figure}
\\
\begin{figure}[tbh]
  \begin{center}
   \subfloat[$\poly{1}$ elements] 
    {\includegraphics[width=.475\linewidth]{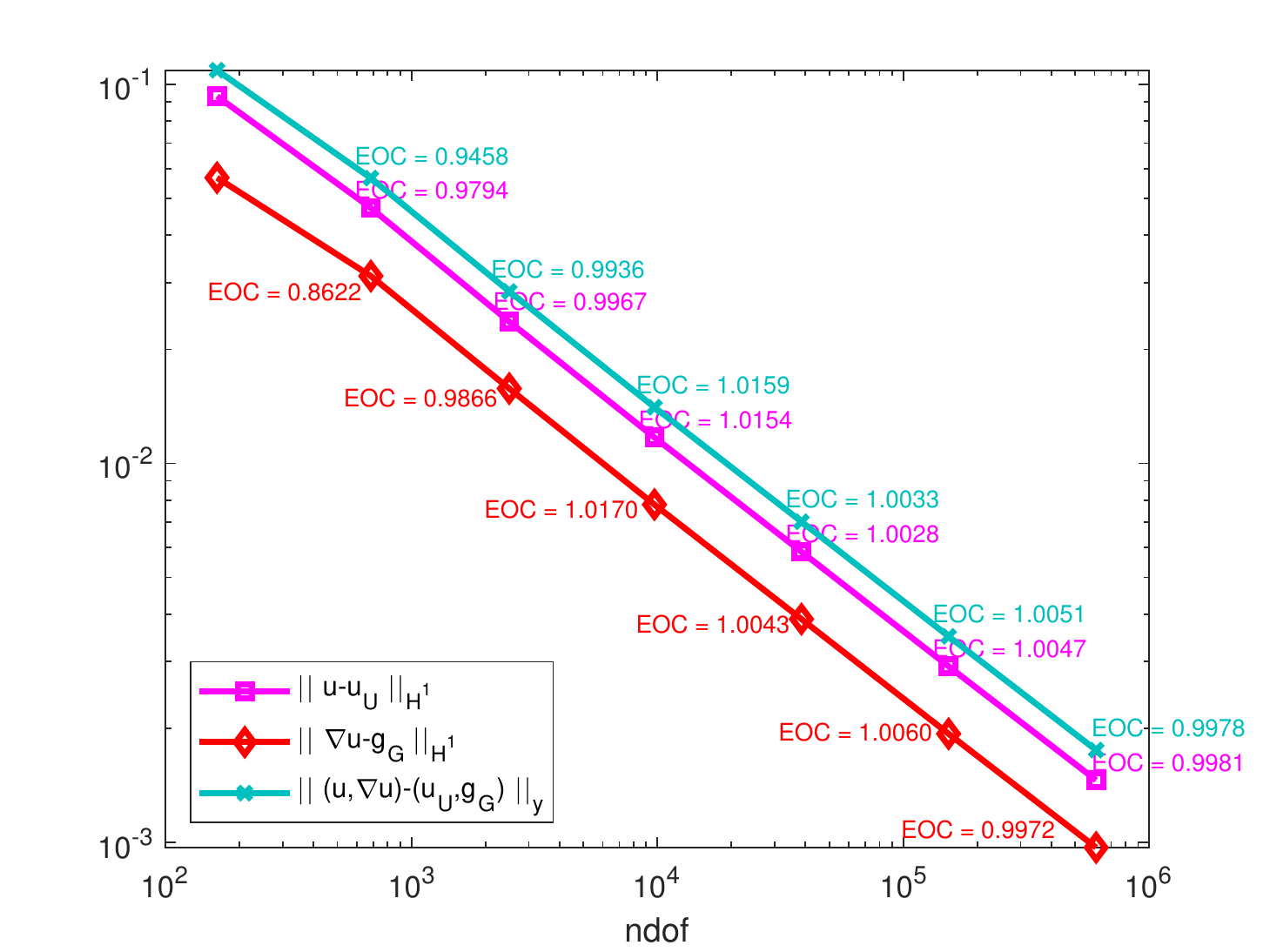}}
    \subfloat[$\poly{2}$ elements] 
    {\includegraphics[width=.475\linewidth]{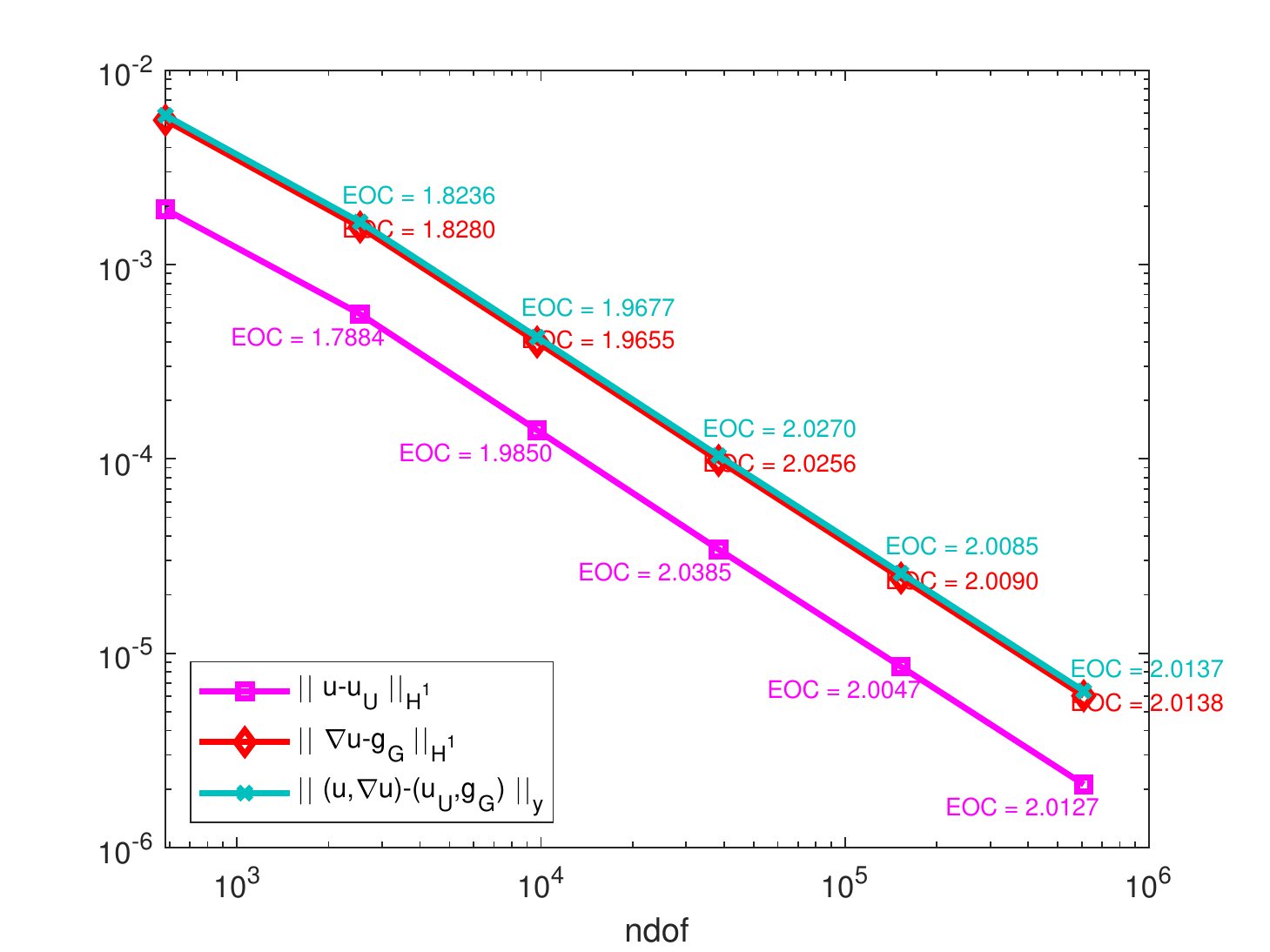}}
    \caption{%
      Experimental order of convergence (EOC) for the Monge--Ampère test problem with $R=2$.}
  \label{fig:MA-R-2}
  \end{center}
\end{figure}
\end{problem}
\begin{problem}[{Hamilton--Jacobi--Bellman test}]
\label{sec:hjb-test-problem}
Consider problem~(\ref{eq:HJB}) and let $\mathcal{A}=[0,2\pi]$,
\begin{equation}
  \label{coef:test-HJB}
  \mat{A}^\alpha(\vec x) = 
  \begin{bmatrix}
    \cos(\alpha) & \sin(\alpha)
    \\
    -\sin(\alpha) & \cos(\alpha)
  \end{bmatrix}
  \begin{bmatrix}
    1+(x_1^2 + x_2^2) & 0.005
    \\
    0.005 & 1.01-(x_1^2 + x_2^2)
  \end{bmatrix}
  \begin{bmatrix}
    \cos(\alpha) & -\sin(\alpha)
    \\
    \sin(\alpha) & \cos(\alpha)
  \end{bmatrix},
\end{equation}
\begin{equation}
\label{data:test-HJB}
\vec b^\alpha =0, 
\quad
c^\alpha = 2 - 0.5(\cos(2\alpha) + \sin(2\alpha)),
\quad 
f^\alpha = \linop L^\alpha u  -(1- \cos(2\alpha - \pi(x_1+x_2))),
\end{equation}
with the exact solution 
\begin{equation}
\label{fanc:exact-solution-HJB}
u(\vec x)
  :=
  \begin{cases}
    r(\vec x)^{5/3} 
    (1-r(\vec x))^{5/2}
    \sin(\varphi(\vec x))^{5/2}
    &\text{ if } 0<r(\vec x) \leq 1 \text{ and }, 0<  \varphi(\vec x) < 3\pi /2,
    \\
    0
    &\text{ otherwise, }
  \end{cases}
\end{equation}
$(r(\vec x), \varphi(\vec x))$ are polar coordinates centered in the
origin.  One can check that the near degenerate diffusion
$\mat{A}^\alpha$ together with $\vec b^\alpha$ and $c^\alpha$ satisfy the
Cordes condition (\ref{eqn:def:general-Cordes-condition}) with $\lambda =
1$ and $\varepsilon = 0.0032$. Note that $u \in \sobh{s}$ for any $s
<8/3$. As $u \in \sobh2(\W)$, we do not expect the advantage of the
adaptive scheme over than the uniform refinement for $\sobh1(\W)$-norm
of the error of $\fespacefun u u$; it is shown in Figure
\ref{fig:HJB-aniso-circle-quad-u}.  But since $\grad u$ does
not have such smoothness, we observe the superiority of the adaptive
scheme for $\sobh1(\W)$-norm of the error of $\vecfespacefun g g$ (and
$\linspace{Y}$-norm of the error of $(\fespacefun u u,
\vecfespacefun g g)$) in Figure \ref{fig:HJB-aniso-circle-quad-g} (and
\ref{fig:HJB-aniso-circle-quad-Y}).
 \begin{figure}[h]
  \begin{center}
    \subfloat[\label{fig:ad-HJB-aniso-circle-6-quad} Last mesh generated by the adaptive algorithm.]{%
      \includegraphics[width=.475\linewidth,trim=20 20 20 20,clip]{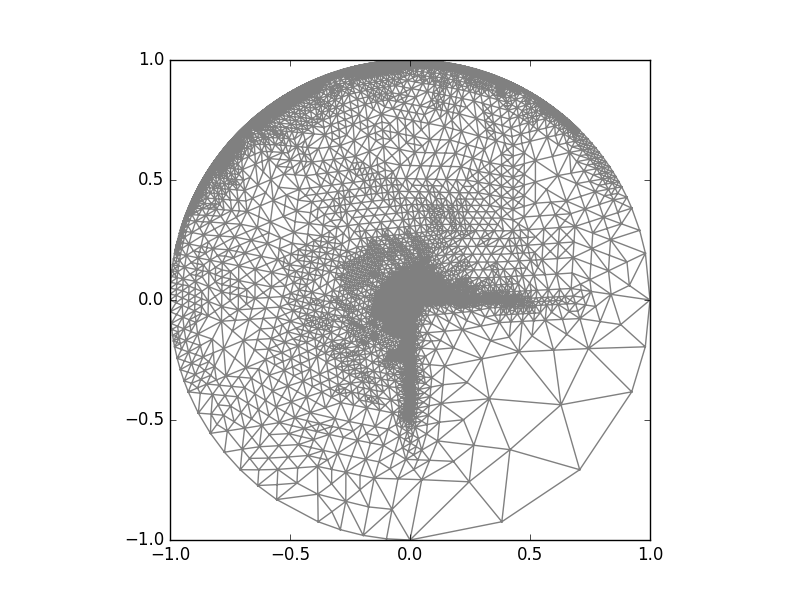}}%
    \subfloat[\label{fig:HJB-aniso-circle-quad-u} Error $\Normon{u-\fespacefun u u}{\sobh1(\W)}$
      convergence of uniform and adaptive methods.]{%
      \includegraphics[width=.475\linewidth]{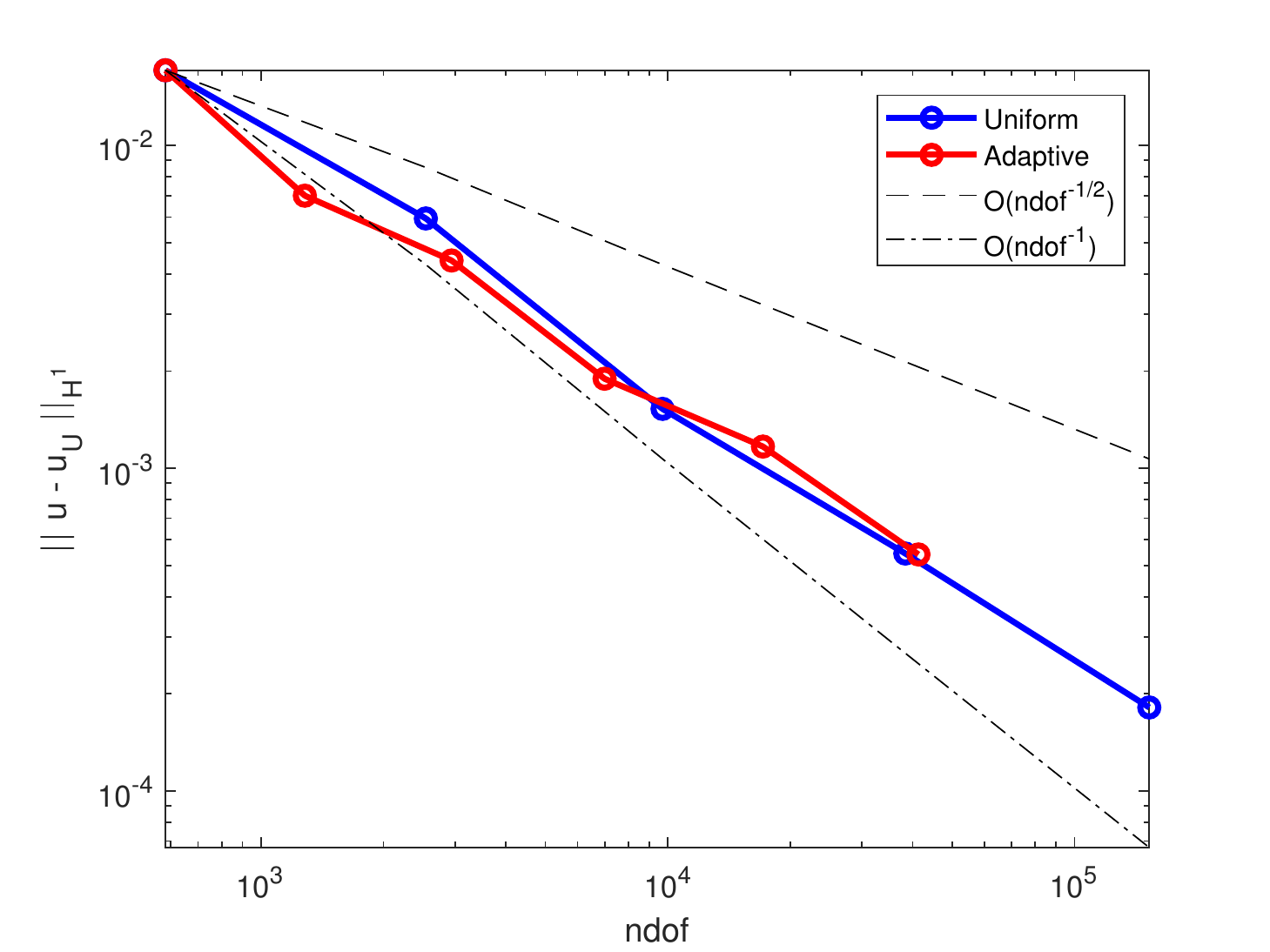}}%
    \\
    \subfloat[\label{fig:HJB-aniso-circle-quad-g}]{%
      \includegraphics[width=.475\linewidth]{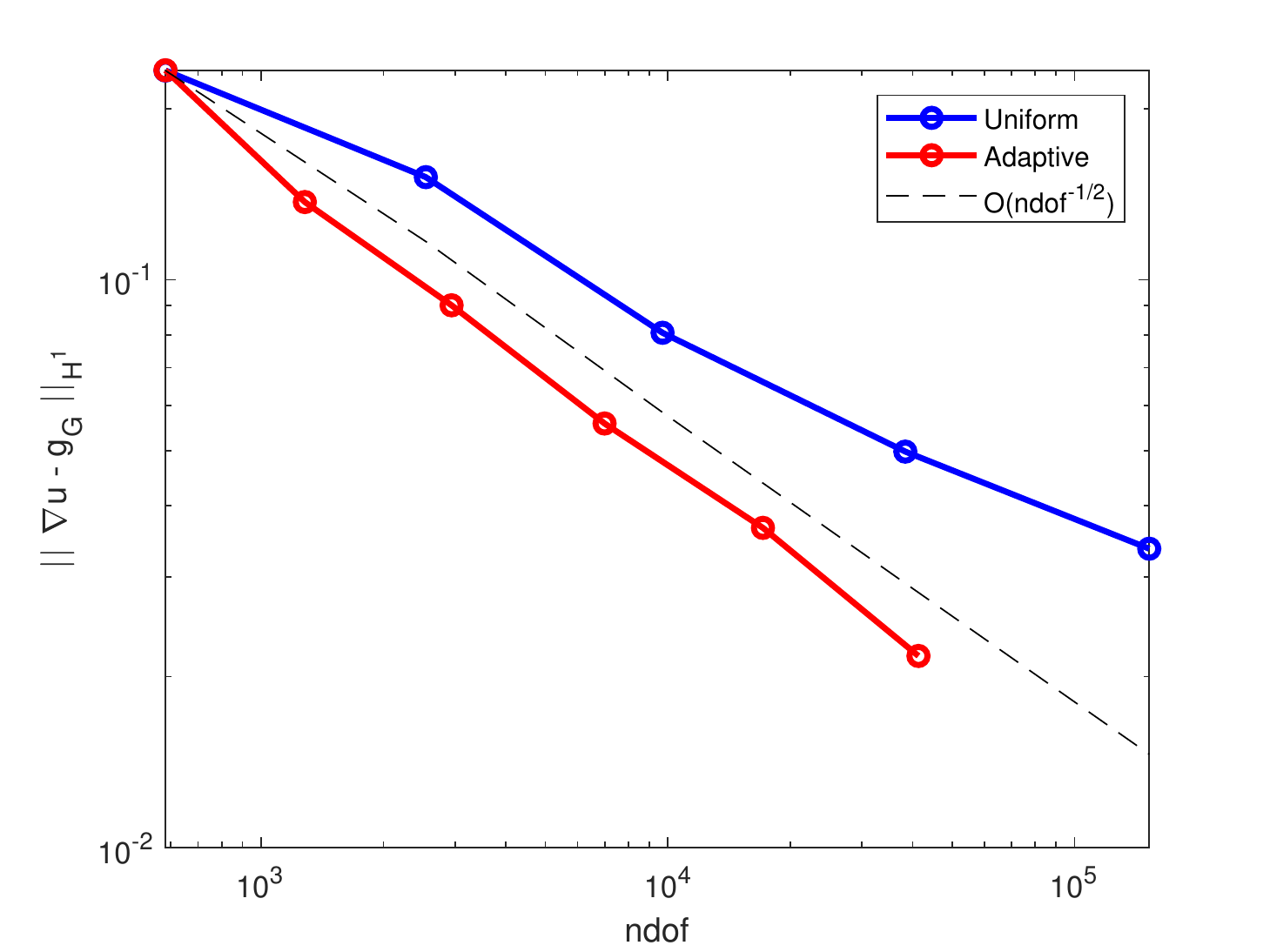}}%
    \subfloat[\label{fig:HJB-aniso-circle-quad-Y}]{%
      \includegraphics[width=.475\linewidth]{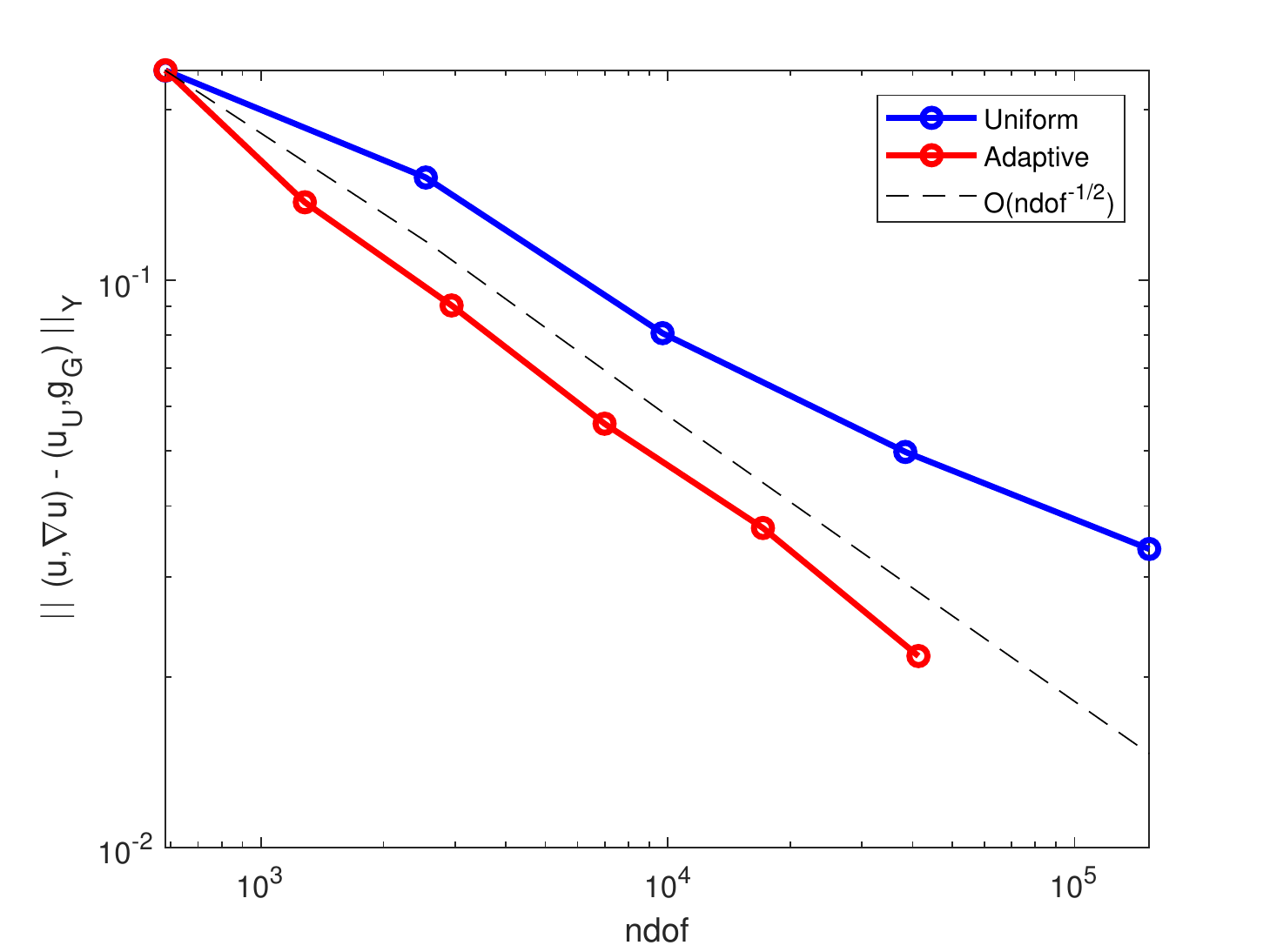}}%
    \caption{%
      \label{fig:HJB}
      Mesh in \ref{fig:ad-HJB-aniso-circle-6-quad} and
      \ref{fig:HJB-aniso-circle-quad-u}--\ref{fig:HJB-aniso-circle-quad-Y}
      show the convergence rate in both the uniform and adaptive
      refinement for the HJB test problem \S\ref{sec:hjb-test-problem}
      with $\poly{2}$ elements.  While the adaptive scheme does not yield any
      noticeable gain for the function value approximation
      ($\Normon{u-\fespacefun u u}{\sobh1(\W)}$), it does so in the
      reconstructed gradient ($\Normonspace{\grad u-\vecfespacefun g g}y$).}
  \end{center}
 \end{figure}
\end{problem}
 \ifthenelse{\boolean{usebibtex}}{
   \bibliographystyle{plainnat}%
 }{
   \renewcommand{\bibname}{References} 
   
 }
\end{document}